\documentclass[a4paper,12pt]{amsart} 
\title[Exact counts in characteristics $2$ and $3$]
{Exact counts of elliptic curves of bounded \\ height over $\mathbb F_q(t)$ in characteristics $2$ and $3$}
\author{Jun--Yong Park}
\date{}

\usepackage[charter]{mathdesign}
\SetMathAlphabet{\mathcal}{normal}{OMS}{cmsy}{m}{n}

\usepackage{amsmath,graphicx,amsthm,amsfonts,verbatim,mathtools,thmtools,tikz-cd,url,tensor}
\usepackage[a4paper,margin=3.3cm, headsep=15pt, headheight=2cm]{geometry}
\usepackage[normalem]{ulem}
\usepackage[shortlabels]{enumitem}
\usepackage[all,cmtip]{xy}
\usepackage{hyperref}
\usepackage{subfiles} 

\newtheorem{thm}{Theorem}[section]

\newtheorem{lem}[thm]{Lemma}
\newtheorem{cor}[thm]{Corollary}
\newtheorem{prop}[thm]{Proposition}
\theoremstyle{definition}

\newtheorem{rmk}[thm]{Remark}

\newcommand{\Zb}{\mathbb{Z}}

\newcommand{\cU}{\mathcal{U}}

\newcommand{\cW}{\mathcal{W}}

\newcommand{\cM}{\mathcal{M}}

\newcommand{\Fb}{\mathbb{F}}

\newcommand{\sm}{\mathrm{sm}}

\DeclareMathOperator{\height}{ht}

\newcommand{\Ac}{\mathcal{A}}

\newcommand{\cE}{\mathcal E}
\newcommand{\Nc}{\mathcal{N}}

\newcommand{\Pb}{\mathbb{P}}

\newcommand{\Gb}{\mathbb{G}}

\newcommand{\Rb}{\mathbb{R}}

\DeclareMathOperator{\Spec}{Spec}

\newcommand{\Aut}{\operatorname{Aut}}
\newcommand{\Sym}{\operatorname{Sym}}

\renewcommand{\setminus}{\smallsetminus}

\begin{document}

    \begin{abstract}
    Let $p\in\{2,3\}$, let $q=p^r$ with $r\geq1$, and put $K=\mathbb F_q(t)$.  We determine the exact weighted and unweighted counts of $K$-isomorphism classes of elliptic curves of bounded Faltings height, equivalently of bounded minimal-discriminant degree.  Writing $B$ for the discriminant height bound, the leading term in each count is of order $B^{5/6}$ and has the same coefficient in every characteristic, whereas the lower-order terms depend substantially on $p$ and on the arithmetic of the constant field.  In the unweighted count, characteristic two produces a term of order $B^{5/12}$ and, when $r$ is even, a term of order $B^{1/4}$, neither of which occurs for $p>3$.  Some characteristic-specific lower-order coefficients are negative.

    These terms reflect two small-characteristic phenomena.  First, the nonsmooth locus of generalized Weierstrass equations contains quasi-elliptic-type strata that are not detected by the rational-singular-section argument in de Jong's count: the unique geometric singular point of the generic fiber may be defined only after a nontrivial purely inseparable extension of $K$.  Second, on the $j=0$ locus, the geometric origin-preserving automorphism groups are nonabelian and contain wild elements, while twisting affects which automorphisms descend to $K$.  Consequently, the passage from weighted to unweighted counts requires a marked-inertia calculation involving conjugacy classes and centralizer weights.

    Our formulas show that nonsmooth-locus corrections, extra-automorphism loci, and the removal of minimality defects collectively account for all lower-order terms.  Together with the characteristic-greater-than-three formulas of Bejleri--Park--Satriano, this completes the exact weighted and unweighted bounded-height enumerations over $\mathbb F_q(t)$ in every characteristic.
    \end{abstract}

    \maketitle
    \vspace{-4ex}


    \section{Introduction}\label{sec:introduction}

Let $p\in\{2,3\}$, let $q=p^r$ with $r\geq1$, and put $K=\Fb_q(t)$.  For every height bound, we determine the exact weighted and unweighted counts of elliptic curves over $K$ up to $K$-isomorphism.  The corresponding formulas in characteristic greater than $3$ were obtained in \cite[Theorem~9.7]{BPS}.  Characteristics $2$ and $3$ are not formal specializations of that calculation: the nonsmooth locus contains additional equations for which the unique geometric singular point is defined only over a nontrivial purely inseparable extension of $K$, while on the $j=0$ locus the geometric origin-preserving automorphism groups are nonabelian and contain wild elements.  Moreover, twisting affects which of these automorphisms descend to $K$.

\medskip

For an elliptic curve $E/K$, let $f\colon X\to\Pb^1_{\Fb_q}$ be its relatively-minimal proper regular model, and let $\sigma\colon\Pb^1_{\Fb_q}\to X$ be the section extending the identity point $O$.  Write $\omega_{X/\Pb^1}$ for the relative dualizing sheaf.  The fundamental line bundle of $f$ and the Faltings height of $E$ are
\[
\mathcal L_f=(R^1f_*\mathcal O_X)^\vee\simeq f_*\omega_{X/\Pb^1},
\qquad
\height(E)=\deg\mathcal L_f,
\]
where the height is a nonnegative integer and $\deg\Delta_{\min}(E)=12\height(E)$.  Throughout, isomorphisms and automorphisms of elliptic curves preserve $O$, while isomorphisms of elliptic surfaces lie over the identity of $\Pb^1_{\Fb_q}$ and preserve the zero section.  We call $E$ \emph{constant} if it is $K$-isomorphic to the base change of an elliptic curve $E_0/\Fb_q$, and \emph{isotrivial} if $E_L\simeq E_0\times_{\Fb_q}L$ for some finite separable extension $L/K$ and some elliptic curve $E_0/\Fb_q$.

\medskip

A generalized Weierstrass equation has \emph{coefficient degree} $d$ if it is written with respect to the line bundle $\mathcal O_{\Pb^1}(d)$, so that $a_i\in H^0(\Pb^1,\mathcal O(id))$ for $i\in\{1,2,3,4,6\}$.  By the \emph{resolved model} of a generalized Weierstrass equation we mean the minimal desingularization of its associated Weierstrass surface; this exists in every characteristic by resolution of singularities for excellent surfaces \cite{Lipman}.  The resolved model need not be relatively-minimal.  If the generic fiber is smooth of exact Faltings height $n$, de Jong's construction compares the resolved model, equipped with its lifted zero section, with the relatively-minimal model and produces an effective blowup divisor $D$ on $\Pb^1$ such that
\[
d=n+\deg D,
\qquad
e=d-n=\deg D.
\]
By analogy with \cite[Definition~4.31]{BPS}, we call $e$ the \emph{minimality defect}.  The divisor $D$ plays the role of the effective quotient divisor in \cite[Lemma~4.32]{BPS}.

\medskip

The starting point is de Jong's weighted-counting framework, in which the class $[E/K]$ has weight $1/|\Aut_K(E,O)|$.  The rational-singular-section parametrization in \cite[Section~4.11]{dJ2} contributes $q^{8d+4}$ coefficient-degree-$d$ equations whose generic fiber has a $K$-rational singular point.  Over the imperfect field $K$, this does not exhaust the nonsmooth locus in characteristics $2$ and $3$.  For certain equations, the unique geometric singular point is defined only after a nontrivial purely inseparable extension of $K$.  The normalization of its closure maps purely inseparably to $\Pb^1$ with degree $p$, rather than giving a section.  Proposition~\ref{prop:nonsmooth} determines the additional contribution, and Theorem~\ref{thm:weighted} gives the corrected weighted count.

\medskip

To obtain the unweighted count, we must also account for $K$-defined automorphisms.  Every elliptic curve has the two distinct central automorphisms $1$ and $[-1]$.  Over an algebraic closure, automorphisms beyond these occur only when $j=0$; the origin-preserving automorphism group then has order $12$ in characteristic three and order $24$ in characteristic two \cite[Appendix~A, Proposition~1.2(c)]{Silverman}.  As shown in Remark~\ref{rmk:etale-aut}, the corresponding automorphism group schemes are nevertheless finite and \'etale.

\medskip

We therefore count marked pairs $(E,\alpha)$ with $\alpha\in\Aut_K(E,O)\setminus\{1,[-1]\}$.  The automorphism group of $(E,\alpha)$ is the centralizer of $\alpha$ in $\Aut_K(E,O)$, so its groupoid weight is
\[
\frac{1}{|C_{\Aut_K(E,O)}(\alpha)|}.
\]
Within each characteristic, the curves on the $j=0$ locus have isomorphic geometric automorphism groups, but twisting changes which automorphisms descend to $K$.

\medskip

At the level of normalized equations, the admissible coordinate changes preserving the normalized form constitute the \emph{residual coordinate group}.  A \emph{fixed pair} consists of a normalized coefficient tuple together with an element of this group that fixes it.  The characteristic-three order-four loci and the characteristic-two type-(II) loci are individually defined over $\Fb_q$ exactly when $r$ is even.  These calculations initially count coefficient-degree-$d$ equations.  We retain the marked automorphism while passing to the relatively-minimal model and then invert the convolution over effective $\Fb_q$-rational divisors $D$.  The intrinsic divisor is essential because divisibility conditions imposed on a chosen normalized coefficient tuple are not preserved by the residual coordinate translations in characteristics $2$ and $3$.

\medskip

For $q=p^r$, set
\begin{equation}\label{eq:parity-indicator}
\epsilon_r=
\begin{cases}
1,&r\text{ even},\\
0,&r\text{ odd}.
\end{cases}
\end{equation}
For every integer $m\geq2$, put
\begin{equation}\label{eq:alpha-beta}
\alpha_m(q)=\sum_{i=0}^{m-2}q^{-i}=\frac{q^{m-1}-1}{q^{m-2}(q-1)},
\qquad
\beta_m(q)=\sum_{i=1}^{m-1}q^{-i}=\frac{q^{m-1}-1}{q^{m-1}(q-1)}=\frac{\alpha_m(q)}{q}.
\end{equation}

\medskip

Define the multiplicative discriminant height by
\[
H_\Delta(E)=q^{\deg\Delta_{\min}(E)}=q^{12\height(E)}.
\]
For $B\geq1$, define the weighted and unweighted counting functions by
\[
\Nc^w(K,B)
=
\sum_{\substack{[E/K]\\ H_\Delta(E)\leq B}}
\frac{1}{|\Aut_K(E,O)|},
\qquad
\Nc(K,B)
=
\#\{[E/K]:H_\Delta(E)\leq B\}.
\]
Here $[E/K]$ ranges over $K$-isomorphism classes of elliptic curves, and all isomorphisms and automorphisms preserve $O$.

\begin{thm}\label{thm:main}
Let $p\in\{2,3\}$ and $r\geq1$, and put $q=p^r$ and $K=\Fb_q(t)$.  For $M\in\Zb_{\geq1}$, put $B=q^{12M}$.

\medskip

If $p=3$, then
\begin{align*}
\Nc^w(K,B)
={}&
\frac{q^9-1}{q^8-q^7}B^{5/6}
-
\beta_4(q)B^{1/3}
-
B^{1/6}
-
1,\\
\Nc(K,B)
={}&
2\frac{q^9-1}{q^8-q^7}B^{5/6}
+
2\alpha_6(q)B^{1/2}\\
&+
2\bigl(
\epsilon_r\alpha_4(q)-2\beta_4(q)
\bigr)B^{1/3}
-
2B^{1/6}
+
2\epsilon_r
-
2.
\end{align*}

If $p=2$, then
\begin{align*}
\Nc^w(K,B)
={}&
\frac{q^9-1}{q^8-q^7}B^{5/6}
-
\beta_6(q)B^{1/2}
-
B^{1/6},\\
\Nc(K,B)
={}&
2\frac{q^9-1}{q^8-q^7}B^{5/6}
+
2\bigl(
2\epsilon_r\alpha_6(q)-\beta_6(q)
\bigr)B^{1/2}\\
&+
\alpha_5(q)B^{5/12}
-
\beta_4(q)B^{1/3}
-
4\epsilon_r\beta_3(q)B^{1/4}
-
2B^{1/6}.
\end{align*}

\medskip

At height zero, equivalently $B=1$,
\[
\Nc^w(K,1)=q,
\qquad
\Nc(K,1)=
\begin{cases}
2q+2+2\epsilon_r,&p=3,\\
2q+1+4\epsilon_r,&p=2.
\end{cases}
\]

For arbitrary $B\in\Rb_{\geq1}$, let $B^\ast$ be the largest element of $q^{12\Zb_{\geq0}}$ not exceeding $B$.  Then $\Nc^w(K,B)=\Nc^w(K,B^\ast)$ and $\Nc(K,B)=\Nc(K,B^\ast)$, so the preceding formulas give the exact weighted and unweighted counts for every $B\in\Rb_{\geq1}$.
\end{thm}

\smallskip

For $M\geq1$ and $B=q^{12M}$, the leading $B^{5/6}$-term of the unweighted count $\Nc(K,B)$ is twice the corresponding term of the weighted count $\Nc^w(K,B)$.  Across the weighted and unweighted formulas, the lower-order terms have three sources: the corrected subtraction of nonsmooth equations, the marked extra-automorphism loci, and the removal of the minimality defect.  The factor $\epsilon_r$ records whether the relevant geometric automorphism loci descend individually to $\Fb_q$.  Section~\ref{sec:twists} gives the full term-by-term explanation.

\medskip

At height zero, the $j=0$ locus consists of $4+2\epsilon_r$ $K$-isomorphism classes represented by constant curves in characteristic three and $3+4\epsilon_r$ in characteristic two, in agreement with \cite[Proposition~2.2]{KST} and \cite[Proposition~3.2]{KST}, respectively.

\medskip

The proof uses generalized Weierstrass equations rather than the isomorphism $\overline{\mathcal M}_{1,1}\simeq\mathcal P(4,6)$ over $\Zb[1/6]$ \cite[Section~1]{BPS}, where $\mathcal P(4,6)=[(\mathbb A^2\setminus\{0\})/\mathbb G_m]$ for the action of weights $4$ and $6$.  This isomorphism does not extend to characteristics $2$ or $3$.  In characteristic $2$, the generic stabilizer of the corresponding scalar action is the nonreduced group scheme $\mu_2$, whereas inversion in the origin-preserving automorphism group generates a constant \'etale subgroup of order two.  In characteristic $3$, the wild order-three automorphisms arise from residual translations of the $x$-coordinate rather than from root-of-unity stabilizers of the scalar action on $\mathcal P(4,6)$; composing these translations with inversion gives the order-six automorphisms.  The proof instead uses the coordinate-change groupoids of generalized Weierstrass equations and their resolved models, while de Jong's invariant blowup divisor passes from coefficient degree to exact Faltings height.

\medskip

The nonsmooth correction is logically distinct from the calculation of wild automorphisms.  Wildness concerns automorphisms whose orders are divisible by the characteristic, whereas the imperfection of $K$ allows the unique geometric singular point to be defined only over a nontrivial purely inseparable extension of $K$.

\subsection{Outline of the paper}

Section~\ref{sec:weighted} corrects the nonsmooth equation count, derives the weighted count, and formulates the marked-inertia correction.  Sections~\ref{sec:char3-fixed} and \ref{sec:char2-fixed} normalize the $j=0$ equations and classify the extra fixed pairs in characteristics $3$ and $2$.  Section~\ref{sec:defect} transfers these ambient marked counts from coefficient degree to exact Faltings height using the intrinsic blowup divisor.  Section~\ref{sec:count} combines the weighted and marked-inertia counts and proves Theorem~\ref{thm:main}.  Finally, Section~\ref{sec:twists} compares the height-zero terms with the finite-field twist classification and explains the geometric origins of the lower-order terms.

\section{Weighted counts and marked inertia}\label{sec:weighted}

Call a groupoid \emph{essentially finite} if it has finitely many isomorphism classes and every object has a finite automorphism group.  For any groupoid $\mathcal G$, write $\pi_0(\mathcal G)$ for its set of isomorphism classes.  If $\mathcal G$ is essentially finite, define its weighted cardinality by
\[
\#_q\mathcal G
=
\sum_{[\xi]\in\pi_0(\mathcal G)}
\frac{1}{|\Aut(\xi)|}.
\]

For $n\geq0$, let $\cE_n(q)$ be the groupoid of elliptic curves $(E,O)/K$ of Faltings height $n$, with morphisms given by origin-preserving $K$-isomorphisms.  After choosing $\mathcal L_f\simeq\mathcal O_{\Pb^1}(n)$, every object of $\cE_n(q)$ admits a generalized Weierstrass equation with $a_i\in H^0(\Pb^1,\mathcal O(in))$ for $i\in\{1,2,3,4,6\}$.  Since these coefficient spaces are finite, only finitely many isomorphism classes occur, while Remark~\ref{rmk:etale-aut} gives finiteness of the automorphism groups. Thus $\cE_n(q)$ is essentially finite.

Define $I_n(q)=\#_q\cE_n(q)$ and $U_n(q)=|\pi_0(\cE_n(q))|$.  Thus $I_n(q)$ is the weighted count, while $U_n(q)$ is the number of $K$-isomorphism classes.

Let $\cE_n^{\mathrm{ex}}(q)$ be the groupoid whose objects are pairs $(E,\alpha)$, where $(E,O)$ is an object of $\cE_n(q)$ and $\alpha\in\Aut_K(E,O)\setminus\{1,[-1]\}$.  A morphism $\phi\colon(E,\alpha)\to(E',\alpha')$ is an origin-preserving $K$-isomorphism satisfying $\phi\alpha\phi^{-1}=\alpha'$.  This groupoid is again essentially finite: there are finitely many such pairs up to isomorphism, and the automorphism group of $(E,\alpha)$ is a subgroup of the finite group $\Aut_K(E,O)$.  Put $X_n(q)=\#_q\cE_n^{\mathrm{ex}}(q)$.

Allowing $\alpha$ to range over all of $\Aut_K(E,O)$ gives the full inertia groupoid of $\cE_n(q)$.  Thus $\cE_n^{\mathrm{ex}}(q)$ is the full subgroupoid obtained by excluding the two markings $1$ and $[-1]$ over every object.  These markings are distinct even in characteristic two, as explained in Remark~\ref{rmk:etale-aut}.  Moreover, every origin-preserving $K$-isomorphism $\phi\colon(E,O)\to(E',O')$ is a group isomorphism \cite[Chapter~III, Theorem~4.8]{Silverman}, and hence $\phi\circ[-1]_E=[-1]_{E'}\circ\phi$.  Consequently, $[-1]$ is central in $\Aut_K(E,O)$, and conjugation by an isomorphism preserves the identity, inversion, and their complement.

The following lemma is the groupoid-cardinality form of orbit--stabilizer and, for quotient groupoids, the fixed-pair mechanism underlying the Cauchy--Frobenius lemma.  Equivalently, it expresses the identity between the weighted point count of an inertia stack and the unweighted point count of the original stack \cite[Theorem~1.1]{HP2}.  We include a direct proof to separate the identity, inversion, and extra-automorphism contributions and to record their centralizer weights.

\begin{lem}\label{lem:inertia}
For every $n\geq0$,
\[
U_n(q)=2I_n(q)+X_n(q).
\]
\end{lem}

\begin{proof}
Fix an isomorphism class $[E]$ and put $A_E=\Aut_K(E,O)$.  The groupoid of pairs $(E',\alpha)$ with $E'\simeq E$ and $\alpha\in\Aut_K(E',O')$ is equivalent to the conjugation action groupoid $A_E\mathbin{//}A_E$.  Indeed, an isomorphism $\psi\colon(E',O')\xrightarrow{\sim}(E,O)$ transports the marking to $\psi\alpha\psi^{-1}\in A_E$, while replacing $\psi$ by $g\psi$ conjugates this element by $g\in A_E$.  Thus the first copy of $A_E$ parametrizes the marked automorphism and the second acts by conjugation.

It follows that the isomorphism classes of $A_E\mathbin{//}A_E$ are the conjugacy classes in $A_E$.  The automorphism group of the object $\alpha\in A_E$ is its centralizer
\[
C_{A_E}(\alpha)=\{g\in A_E:g\alpha=\alpha g\}.
\]
By orbit--stabilizer,
$|[\alpha]_{\mathrm{conj}}|=|A_E|/|C_{A_E}(\alpha)|$.  Consequently,
\[
\begin{aligned}
\#_q\bigl(A_E\mathbin{//}A_E\bigr)
&=
\sum_{[\alpha]_{\mathrm{conj}}\subset A_E}
\frac{1}{|C_{A_E}(\alpha)|}\\
&=
\frac{1}{|A_E|}
\sum_{[\alpha]_{\mathrm{conj}}\subset A_E}
|[\alpha]_{\mathrm{conj}}|
=1,
\end{aligned}
\]
because the conjugacy classes partition $A_E$.

The identity and inversion are distinct central elements of $A_E$, so their two conjugacy classes contribute $2/|A_E|$.  The remaining classes contribute $1-2/|A_E|=(|A_E|-2)/|A_E|$.  They are precisely the isomorphism classes in the fiber of $\cE_n^{\mathrm{ex}}(q)$ over $[E]$, and the automorphism group of $(E,\alpha)$ is $C_{A_E}(\alpha)$.  Summing this decomposition over $[E]\in\pi_0(\cE_n(q))$ gives the result.
\end{proof}

\begin{prop}
\label{prop:rational-points-of-quotients}
Let $G$ be an algebraic group over $\Fb_q$, and let $X$ be a left $G$-scheme of finite type over $\Fb_q$.  Write $H^1(\Fb_q,G):=H^1_{\mathrm{fppf}}(\Fb_q,G)$ for the pointed set of isomorphism classes of $G$-torsors, and suppose that $H^1(\Fb_q,G)=\{1\}$.  Then the canonical functor $X(\Fb_q)\mathbin{//}G(\Fb_q)\to[X/G](\Fb_q)$ is an equivalence.

Let $Z\subset X\times G$ be the scheme-theoretic fixed-pair incidence defined by $g\cdot x=x$, with $G$ acting by $h\cdot(x,g)=(h\cdot x,hgh^{-1})$.  Then $I([X/G])\simeq[Z/G]$, and
\[
\begin{aligned}
\left|\pi_0\bigl([X/G](\Fb_q)\bigr)\right|
&=
\left|X(\Fb_q)/G(\Fb_q)\right|\\
&=
\frac{1}{|G(\Fb_q)|}
\sum_{g\in G(\Fb_q)}|X(\Fb_q)^g|\\
&=
\frac{|Z(\Fb_q)|}{|G(\Fb_q)|}
=
\#_q I([X/G]).
\end{aligned}
\]
Here $X(\Fb_q)^g$ is the fixed-point set of $g$, and $\#_q I([X/G])$ denotes the weighted cardinality of $I([X/G])(\Fb_q)$.  For every $(x,g)\in Z(\Fb_q)$, its stabilizer in $G(\Fb_q)$ is
\[
C_{\operatorname{Stab}_{G(\Fb_q)}(x)}(g).
\]
\end{prop}

\begin{proof}
The canonical functor
\[
X(\Fb_q)\mathbin{//}G(\Fb_q)\longrightarrow[X/G](\Fb_q)
\]
is always fully faithful, and its essential image consists precisely of the objects whose underlying $G$-torsor is trivial.  Since $H^1(\Fb_q,G)=\{1\}$, every $G$-torsor is trivial, so this functor is essentially surjective and hence an equivalence.  The same argument, with $X$ replaced by $Z$, gives
\[
[Z/G](\Fb_q)\simeq Z(\Fb_q)\mathbin{//}G(\Fb_q).
\]
The inertia presentation $I([X/G])\simeq[Z/G]$ is established in \cite[Corollary~2.5]{HP2}.  The displayed equalities now follow from the Cauchy--Frobenius formula (Burnside's lemma), the identity $|Z(\Fb_q)|=\sum_{g\in G(\Fb_q)}|X(\Fb_q)^g|$, and the groupoid-cardinality formula for $Z(\Fb_q)\mathbin{//}G(\Fb_q)$.  Finally, $h\in G(\Fb_q)$ fixes $(x,g)$ precisely when $h\cdot x=x$ and $hgh^{-1}=g$, equivalently when $h\in C_{\operatorname{Stab}_{G(\Fb_q)}(x)}(g)$.  Thus no commutativity hypothesis on the stabilizer or its centralizer is required. 
\end{proof}

\begin{rmk}
\label{rmk:etale-aut}
Let $(E,O)/K$ be an elliptic curve and put $V=H^0(E,\mathcal O_E(3O))^\vee$.  The complete linear system $|3O|$ embeds $E$ as a plane cubic in $\Pb(V)$.  Under this embedding, the automorphism functor is the closed stabilizer subgroup scheme
\[
G=\underline{\Aut}_K(E,O)
\hookrightarrow
\operatorname{PGL}(V)
\]
of the embedded cubic together with the point $O$.  Hence $G$ is a group scheme of finite type.  For every $K$-algebra $R$, one has $G(R)=\Aut_R(E_R,O_R)$, where $(E_R,O_R)=(E,O)\times_K\Spec R$; in particular, $G(K)=\Aut_K(E,O)$.

The Lie algebra of $G$ is its tangent space at the identity:
\[
\operatorname{Lie}(G)
=T_1G
=
\ker\left(
G\bigl(K[\varepsilon]/(\varepsilon^2)\bigr)
\longrightarrow
G(K)
\right).
\]
It parametrizes first-order automorphisms that reduce to the identity when $\varepsilon=0$.  First-order automorphisms of $E$ correspond to global vector fields, and fixing $O$ forces the vector field to vanish there. An invariant differential trivializes the cotangent bundle of $E$, so $T_E\simeq\mathcal O_E$.  Therefore
\[
\operatorname{Lie}(G)
=H^0\bigl(E,T_E(-O)\bigr)
\simeq
H^0\bigl(E,\mathcal O_E(-O)\bigr)
=0,
\]
because $\mathcal O_E(-O)$ has degree $-1$.

It follows that $G$ is unramified at the identity.  Formation of the Lie algebra commutes with extension of scalars, so $\operatorname{Lie}(G_{\overline K})=0$.  Left translation identifies the tangent space at every geometric point with the tangent space at the identity.  Thus $G_{\overline K}$ is unramified, and hence $G$ is unramified over $K$ by descent.  Since every $K$-scheme is flat over $K$ and a finite-type $K$-scheme is of finite presentation, $G$ is \'etale over $K$.  An \'etale $K$-scheme of finite type is a finite disjoint union of spectra of finite separable extensions of $K$.  Consequently, $G$ is finite \'etale.

In characteristic two, the field elements $-1$ and $1$ coincide, but the inversion morphism $[-1]\colon E\to E$ is not the identity.  For every nonzero integer $m$, multiplication by $m$ is a finite isogeny of degree $m^2$; equivalently, its kernel group scheme $E[m]$ has rank $m^2$.  In particular, $\deg[2]=4$ even in characteristic two.  If $[-1]=[1]$ as morphisms, then the group law would give $[2]=[1]+[1]=[1]+[-1]=[0]$, where $[0]$ is the constant morphism with value $O$.  This contradicts the fact that $[2]$ is finite of degree four.

For $j(E)\neq0$, the geometric automorphism group consists of $1$ and $[-1]$ \cite[Appendix~A, Proposition~1.2(c)]{Silverman}.  Since both are defined over $K$, the automorphism group scheme is the constant group scheme with these two elements.  In characteristic two, this is the constant \'etale group scheme $(\Zb/2\Zb)_K$, not the nonreduced group scheme $\mu_2$.  All weighted counts below use the finite group $G(K)=\Aut_K(E,O)$ of $K$-defined automorphisms.
\end{rmk}

We next revisit the geometric setup underlying de Jong's weighted count and refine its nonsmooth-locus calculation in characteristics $2$ and $3$. A possibly nonminimal generalized Weierstrass equation over $\Pb^1$ is written with respect to a line bundle $\mathcal L$, with $a_i\in H^0(\Pb^1,\mathcal L^{\otimes i})$ for $i\in\{1,2,3,4,6\}$.  We call $d=\deg\mathcal L$ its \emph{coefficient degree}.  Since $\mathcal L\simeq\mathcal O_{\Pb^1}(d)$, its coefficients satisfy
\[
a_i\in H^0\bigl(\Pb^1,\mathcal O(id)\bigr),
\qquad i\in\{1,2,3,4,6\}.
\]
The coefficient degree need not equal the Faltings height, since the associated resolved model need not be relatively-minimal.  We use $n$ for the exact Faltings height and reserve $M$ for the cutoff in cumulative counts.

If the relatively-minimal model obtained from a coefficient-degree-$d$ equation has height $n$, then $n\leq d$.  De Jong's effective blowup divisor $D$ satisfies
\[
e:=\deg D=d-n.
\]
By analogy with \cite[Definition~4.31]{BPS}, we call $e$ the \emph{minimality defect}.  The blowups encoded by $D$ do not change the generic elliptic curve; $e$ records the nonminimality of the chosen Weierstrass model over $\Pb^1$.

Since $h^0(\Pb^1,\mathcal O(id))=id+1$, the full coefficient space has $q^{16d+5}$ elements.  De Jong's singular-section argument counts $q^{8d+4}$ equations whose generic fiber has a $K$-rational singular point \cite[Section~4.11]{dJ2}.  Over the imperfect field $K=\Fb_q(t)$, however, the unique geometric singular point may be defined only over a purely inseparable extension, so the rational-singular count need not exhaust the nonsmooth locus.  As the proof below shows, the additional equations can occur only in the locus $a_1=a_3=0$ in characteristic two or, after completing the square, in the locus $a_2=a_4=0$ in characteristic three.  We call these the \emph{quasi-elliptic-type strata}.  After extension to $\overline K$ and an admissible Weierstrass coordinate change, their geometric generic fibers are cuspidal cubics of the form $Y^2=X^3$.

\begin{prop}\label{prop:nonsmooth}
Let $p$ be a prime, let $r\geq1$, put $q=p^r$, and let $d\geq0$.  If $p\in\{2,3\}$, the number of coefficient-degree-$d$ generalized Weierstrass equations over $\Pb^1_{\Fb_q}$ whose generic fiber is not smooth equals
\[
\begin{cases}
q^{8d+4}+q^{10d+3}-q^{6d+3},&p=3,\\[1mm]
q^{8d+4}+q^{12d+3}-q^{7d+3},&p=2.
\end{cases}
\]
If $p\geq5$, the number is $q^{8d+4}$.  In particular, for $d=0$ the answer is $q^4$ in every characteristic.
\end{prop}

\begin{proof}
Write the equation as $F=0$, where
\[
F
=
y^2+a_1xy+a_3y-x^3-a_2x^2-a_4x-a_6,
\qquad
a_i\in H^0\bigl(\Pb^1,\mathcal O(id)\bigr).
\]
The generic fiber is smooth if and only if its discriminant is nonzero in $K$.

We first count the equations whose generic fiber has a $K$-rational singular point.  After homogenization, the generic Weierstrass cubic meets the line at infinity only at $O=[0:1:0]$.  This point is smooth, since the derivative of the homogenized equation with respect to the homogenizing variable equals $1$ at $O$.  Every positive-degree projective component meets the line at infinity, so every geometric component of the cubic passes through $O$.  If the cubic were reducible, at least two components would pass through $O$; if it were nonreduced, a repeated component would pass through $O$.  Either possibility would make $O$ singular.  Thus the cubic is geometrically integral.  By B\'ezout's theorem, a geometrically integral plane cubic has at most one geometric singular point: otherwise, the line joining two distinct singular points would meet the cubic with total intersection multiplicity at least four.

By the valuative criterion of properness over the regular curve $\Pb^1$, a $K$-rational singular point extends uniquely to a section of the associated Weierstrass surface.  The relative singular locus, namely the closed locus where the morphism to $\Pb^1$ is not smooth, contains the generic point of this section and hence the entire section.  The zero section lies in the relative smooth locus, and every fiber meets the line at infinity only in the zero section.  Hence the singular section lies in the affine chart and has the form
\[
(x,y)=(C_2,C_3),
\qquad
C_2\in H^0\bigl(\Pb^1,\mathcal O(2d)\bigr),
\qquad
C_3\in H^0\bigl(\Pb^1,\mathcal O(3d)\bigr).
\]
At $(C_2,C_3)$, the equations $F=F_x=F_y=0$ are equivalent to
\[
\begin{aligned}
a_3&=-2C_3-a_1C_2,\\
a_4&=a_1C_3-3C_2^2-2a_2C_2,\\
a_6&=-C_3^2-a_1C_2C_3+2C_2^3+a_2C_2^2.
\end{aligned}
\]
Since the geometric singular point is unique, $(a_1,a_2,C_2,C_3)$ parametrizes these equations bijectively.  Their number is
\[
q^{(d+1)+(2d+1)+(2d+1)+(3d+1)}=q^{8d+4}.
\]

We use the following elementary observation.  Let $m\geq1$ and let $s$ be a rational section of a line bundle $L$ on $\Pb^1$.  If $s^m\in H^0(\Pb^1,L^{\otimes m})$, then $s$ is global: in a local trivialization at every closed point $x$, one has $m v_x(s)=v_x(s^m)\geq0$.

Suppose first that $p=2$.  On the affine chart, $\partial F/\partial y=a_1x+a_3$.  If $a_1\neq0$, every geometric singular point satisfies
\[
x_0=\frac{a_3}{a_1},
\qquad
a_1y_0=x_0^2+a_4,
\]
and is therefore $K$-rational.  If $a_1=0$ and $a_3\neq0$, then $\partial F/\partial y=a_3$ is nonzero on the generic fiber, so the generic fiber is smooth.  If $a_1=a_3=0$, the singularity equations over $\overline K$ are
\[
x_0^2=a_4,
\qquad
y_0^2=a_2a_4+a_6.
\]
They have a unique solution in $\overline K^2$, so all $q^{(2d+1)+(4d+1)+(6d+1)}=q^{12d+3}$ equations in this stratum are nonsmooth.

Within this stratum, the singular point is $K$-rational precisely when there are rational sections $C_2$ of $\mathcal O(2d)$ and $C_3$ of $\mathcal O(3d)$ such that
\[
a_4=C_2^2,
\qquad
a_6=C_3^2+a_2C_2^2.
\]
The preceding observation shows that both sections are global.  Since the squaring maps on their spaces of sections are injective, the triples $(a_2,C_2,C_3)$ parametrize exactly $q^{(2d+1)+(2d+1)+(3d+1)}=q^{7d+3}$ members of this stratum already counted by the rational-singular parametrization.  Inclusion--exclusion therefore gives
\[
q^{8d+4}+q^{12d+3}-q^{7d+3}.
\]

Now suppose that $p=3$.  Completing the square removes $a_1$ and $a_3$. For each fixed pair $(a_1,a_3)$, this is an affine-linear bijection on the remaining coefficient space, and there are $q^{(d+1)+(3d+1)}=q^{4d+2}$ choices of $(a_1,a_3)$.  Renaming the transformed coefficients gives
\[
y^2=f(x),
\qquad
f=x^3+a_2x^2+a_4x+a_6.
\]
A singular point satisfies
\[
y_0=0,
\qquad
f(x_0)=f'(x_0)=0,
\qquad
f'=2a_2x+a_4.
\]
If $a_2\neq0$, every repeated root is $x_0=-a_4/(2a_2)\in K$.  If $a_2=0$ and $a_4\neq0$, then $f'=a_4$ is nonzero, so the generic fiber is smooth.  If $a_2=a_4=0$, then $f=x^3+a_6$.  For each fixed pair $(a_1,a_3)$, all $q^{6d+1}$ choices of $a_6$ in this stratum are nonsmooth.  The singular point is $K$-rational precisely when $a_6=-C_2^3$ for a rational section $C_2$ of $\mathcal O(2d)$.  The preceding observation shows that $C_2$ is global, and injectivity of the cubing map gives $q^{2d+1}$ such equations.

For a rational singular point with $x=C_2$, the equations $f'(C_2)=f(C_2)=0$ give
\[
a_4=-2a_2C_2,
\qquad
a_6=-C_2^3+a_2C_2^2.
\]
For each fixed $(a_1,a_3)$, the rational-singular locus has $q^{4d+2}$ elements, the stratum $a_2=a_4=0$ has $q^{6d+1}$ elements, and their intersection, obtained by setting $a_2=0$, has $q^{2d+1}$ elements. Restoring the $q^{4d+2}$ choices of $(a_1,a_3)$, inclusion--exclusion gives
\[
q^{4d+2}
\bigl(q^{4d+2}+q^{6d+1}-q^{2d+1}\bigr)
=q^{8d+4}+q^{10d+3}-q^{6d+3}.
\]

For $p\geq5$, completing the square and depressing the cubic gives
\[
y^2=f(x),
\qquad
f=x^3+Ax+B.
\]
These coordinate changes are defined over $K$ and preserve the rationality of singular points.  A singular point has $y=0$ and is a common root of $f$ and $f'=3x^2+A$.  If $A=0$, the only root of $f'$ is $0\in K$.  If $A\neq0$, then $f'$ is separable.  An irrational common root of $f$ and $f'$ would have a distinct conjugate that was also a common root, giving two distinct repeated roots of a cubic, which is impossible. Hence every nonsmooth equation has a $K$-rational singular point, and its number is $q^{8d+4}$.

Finally, when $d=0$, Frobenius is bijective on $\Fb_q$, so the small-characteristic strata add no equations beyond the rational-singular locus.  Both formulas reduce to $q^4$.
\end{proof}

\begin{rmk}
\label{rmk:dejong-gap}
The assertion in \cite[Section~4.11]{dJ2} that the unique geometric singular point extends to a section over $\Pb^1$ requires modification. The count $q^{8d+4}$ obtained there includes precisely the equations whose singular point is $K$-rational.  Over the imperfect field $K=\Fb_q(t)$, the unique geometric singular point may instead be defined over a purely inseparable extension.

For example, $y^2=x^3+t$ has geometric singular point $(0,t^{1/2})$ in characteristic two and $(-t^{1/3},0)$ in characteristic three.  Let $v_0$ denote the order of vanishing at $t=0$.  Since $v_0(t)=1$, whereas $v_0(f^p)=p\,v_0(f)$ for every $f\in K^\times$, one has $t\notin K^p$. Thus, in either characteristic, the singular point is not $K$-rational. In either characteristic, the normalization of its closure maps to $\Pb^1$ by a purely inseparable morphism of degree $p$, rather than giving a section.

Consequently, for $d\geq1$ and $p\in\{2,3\}$, the term $q^{8d+4}$ used for the type-(b) locus in the proof of \cite[Proposition~4.12]{dJ2} must be replaced by the corresponding formula in Proposition~\ref{prop:nonsmooth}. Only the type-(b) nonsmooth term is affected; the type-(a) blowup-divisor construction in \cite[Section~4.11]{dJ2} is unchanged.  The additional contribution comes from the quasi-elliptic-type strata without a rational singular section.  Quasi-elliptic fibrations occur only in characteristics two and three \cite{BMIII}, and the contribution is
\[
\begin{cases}
q^{10d+3}-q^{6d+3},&p=3,\\
q^{12d+3}-q^{7d+3},&p=2.
\end{cases}
\]
For $d\geq1$, this exceeds the original $q^{8d+4}$ rational-singular term, but remains lower order than the $q^{16d+5}$ full coefficient space.  It therefore changes the lower-order terms without changing the leading weighted asymptotic.
\end{rmk}

Let $W_d^{\sm}(q)$ be the finite set of coefficient-degree-$d$ generalized Weierstrass equations whose generic fiber is smooth, and let $G_d(q)$ be the group of global Weierstrass coordinate changes acting on it.  Explicitly, these changes have the form
\[
\begin{gathered}
x=u^2x'+\rho,
\qquad
y=u^3y'+u^2\eta x'+\tau,\\
u\in\Fb_q^\times,
\quad
\rho\in H^0\bigl(\Pb^1,\mathcal O(2d)\bigr),
\quad
\eta\in H^0\bigl(\Pb^1,\mathcal O(d)\bigr),
\quad
\tau\in H^0\bigl(\Pb^1,\mathcal O(3d)\bigr).
\end{gathered}
\]
Define the action groupoid $\cW_d^{\sm}(q)=W_d^{\sm}(q)\mathbin{//}G_d(q)$ and put $R_d(q)=\#_q\cW_d^{\sm}(q)$.  Thus $R_d(q)$ is the weighted count of generically smooth coefficient-degree-$d$ equations modulo global coordinate change, before relative minimality is imposed.  Since this is a finite action groupoid,
\[
R_d(q)=\frac{|W_d^{\sm}(q)|}{|G_d(q)|}.
\]

The choices of $u,\rho,\eta$, and $\tau$ give
\[
|G_d(q)|
=(q-1)q^{(2d+1)+(d+1)+(3d+1)}
=(q-1)q^{6d+3}
\]
\cite[Sections~4.8--4.9]{dJ2}.  Proposition~\ref{prop:nonsmooth} therefore
gives
\begin{equation}\label{eq:rawall}
R_d(q)
=
\begin{cases}
\dfrac{q^{10d+2}-q^{4d}-q^{2d+1}+1}{q-1},
&p=3,\\[3mm]
\dfrac{q^{10d+2}-q^{6d}-q^{2d+1}+q^d}{q-1},
&p=2.
\end{cases}
\end{equation}
In both characteristics, $R_0(q)=q$.

For $e\geq0$, put
\[
c_e(q)
=
\bigl|\Sym^e(\Pb^1)(\Fb_q)\bigr|
=
\frac{q^{e+1}-1}{q-1}.
\]
Here $\Sym^e(\Pb^1)(\Fb_q)$ parametrizes effective $\Fb_q$-rational divisors of degree $e$, possibly supported at nonrational closed points.  The formula follows from $\Sym^e(\Pb^1)\simeq\Pb^e$.

The blowup construction of \cite[Section~4.11]{dJ2} associates to a relatively-minimal pair $(X,\sigma)$ of exact height $n\leq d$ and an effective divisor
\[
D=\sum_{x\in|\Pb^1|}e_x\,x
\in
\Sym^{d-n}(\Pb^1)(\Fb_q)
\]
a coefficient-degree-$d$ model.  Here $(X,\sigma)$ is the relatively-minimal proper regular elliptic surface with zero section associated with an object $(E,O)$ of $\cE_n(q)$, and restriction to the generic fiber identifies $\Aut_{\Pb^1}(X,\sigma)\simeq\Aut_K(E,O)$.

For fixed $(X,\sigma)$, the resulting model is determined up to isomorphism by $D$.  Conversely, every generically smooth coefficient-degree-$d$ model arises from such a pair and a divisor $D$ \cite[Section~4.11]{dJ2}; Lemma~\ref{lem:blowup-weierstrass} shows that successively contracting the vertical $(-1)$-curves recovers both $(X,\sigma)$ and $D$ uniquely and preserves the automorphism group.  Thus, for fixed $(X,\sigma)$, the stabilizer weight is independent of $D$.  Together with the stabilizer-weighted coordinate count of \cite[Lemma~4.9]{dJ2}, this gives
\begin{equation}\label{eq:dejongrecursion}
R_d(q)
=
\sum_{n=0}^{d}c_{d-n}(q)I_n(q).
\end{equation}
This is the finite-field weighted-point realization of the minimality-defect stratification of \cite[Proposition~6.1 and Corollary~6.2]{BPS}; the corresponding motivic generating-series factorization is \cite[Lemma~8.6, equation~(6)]{BPS}.

\begin{thm}\label{thm:weighted}
For every $M\geq1$,
\[
\sum_{n=0}^M I_n(q)
=
\begin{cases}
\dfrac{q^9-1}{q^8-q^7}\,q^{10M}
-q^{2M}
-\beta_4(q)q^{4M}
-1,
&p=3,\\[3mm]
\dfrac{q^9-1}{q^8-q^7}\,q^{10M}
-q^{2M}
-\beta_6(q)q^{6M},
&p=2.
\end{cases}
\]
At height zero, $I_0(q)=q$ in both characteristics.
\end{thm}

\begin{proof}
Put
\[
R(z)=\sum_{d\geq0}R_d(q)z^d,
\qquad
I(z)=\sum_{n\geq0}I_n(q)z^n.
\]
Since $c_e(q)=1+q+\cdots+q^e$, one has
\[
\sum_{e\geq0}c_e(q)z^e
=
\frac{1}{(1-z)(1-qz)}.
\]
Equation~\eqref{eq:dejongrecursion} therefore gives
\[
R(z)=\frac{I(z)}{(1-z)(1-qz)}.
\]

For a formal power series $H(z)=\sum_{m\geq0}h_mz^m$, write $[z^M]H(z)=h_M$.  Since $(1-z)^{-1}=\sum_{j\geq0}z^j$, the coefficient of $z^M$ in $I(z)/(1-z)$ is $\sum_{n=0}^M I_n(q)$.  The coefficient of $z^M$ in $(1-qz)R(z)$ is $R_M(q)-qR_{M-1}(q)$.  Hence, for $M\geq1$,
\[
\sum_{n=0}^M I_n(q)
=
[z^M]\frac{I(z)}{1-z}
=
[z^M](1-qz)R(z)
=
R_M(q)-qR_{M-1}(q).
\]
Substituting \eqref{eq:rawall} gives the two formulas in the statement. The term $-q^{2M}$ comes from rational singular sections.  The additional quasi-elliptic-type corrections are $-\beta_4(q)q^{4M}-1$ in characteristic three and $-\beta_6(q)q^{6M}$ in characteristic two.  Finally, $c_0(q)=1$, so setting $d=0$ in \eqref{eq:dejongrecursion} gives $I_0(q)=R_0(q)=q$.
\end{proof}

With $I_n(q)$ determined, Lemma~\ref{lem:inertia} reduces the unweighted count to computing $X_n(q)$.  Over $\overline K$, an elliptic curve in characteristic two or three has automorphisms beyond $\{1,[-1]\}$ only when $j=0$ \cite[Appendix~A, Proposition~1.2(c)]{Silverman}.  Since a $K$-defined extra automorphism remains extra after extension of scalars, every object of $\cE_n^{\mathrm{ex}}(q)$ lies on the $j=0$ locus.  Sections~\ref{sec:char3-fixed} and \ref{sec:char2-fixed} give normalized quotient descriptions of the corresponding coefficient spaces and count coefficient-degree-$d$ equations marked by an extra $K$-defined automorphism.  Section~\ref{sec:defect} then passes from coefficient degree to exact Faltings height using de Jong's intrinsic blowup divisor, as formalized in Proposition~\ref{prop:marked-defect}.

\section{Extra automorphisms in characteristic 3}
\label{sec:char3-fixed}

Assume throughout that $q=3^r$ with $r\geq1$.  For $i,d\geq0$, write $V_i(d)=H^0(\Pb^1,\mathcal O(id))$.  We use the same notation for the associated affine space and, when appropriate, for its additive vector group.

\subsection{Normalized equations and the residual coordinate group}

In characteristic three, every generically smooth $j=0$ Weierstrass equation of coefficient degree $d$ admits a normalized form
\begin{equation}\label{eq:char3-normal}
E_{A,C}\colon y^2=x^3+Ax+C,
\qquad
0\neq A\in V_4(d),
\qquad
C\in V_6(d)
\end{equation}
\cite[Section~4.13(b2)]{dJ2}.  The condition $A\neq0$ is equivalent to generic smoothness, since the discriminant is a nonzero scalar multiple of $A^3$.

The \emph{residual coordinate group} consists of the global coordinate changes preserving \eqref{eq:char3-normal}.  It is $G_d^{(3)}=\Gb_m\ltimes V_2(d)$, where $u\in\Gb_m$ acts on $V_2(d)$ by multiplication by $u^2$.  Thus
\begin{equation}\label{eq:char3-group-size}
|G_d^{(3)}(\Fb_q)|=(q-1)q^{2d+1},
\end{equation}
and its group law is
\begin{equation}\label{eq:char3-group-law}
(u_1,s_1)(u_2,s_2)
=
(u_1u_2,u_1^2s_2+s_1).
\end{equation}
The corresponding coordinate change is
\begin{equation}\label{eq:char3-change}
(x,y)\longmapsto(u^2x+s,u^3y),
\end{equation}
which induces the left action
\begin{equation}\label{eq:char3-action}
(u,s)\cdot(A,C)
=
\bigl(u^4A,u^6C-s^3-u^4As\bigr)
\end{equation}
on
$\cU_d^{(3),0}=(V_4(d)\setminus\{0\})\times V_6(d)$.

\begin{prop}
\label{prop:char3-normalized-quotient}
Let $\cW_d^{(3),0}(q)$ be the groupoid of generically smooth $j=0$ Weierstrass equations of coefficient degree $d$ over the fixed base $\Pb^1_{\Fb_q}$.  Its morphisms are global Weierstrass coordinate changes over the identity of $\Pb^1_{\Fb_q}$.  Normalization induces an equivalence
\[
\cW_d^{(3),0}(q)
\simeq
[\cU_d^{(3),0}/G_d^{(3)}](\Fb_q).
\]

Let $\overline X_{A,C}\to\Pb^1$ be the corresponding Weierstrass surface, and let $\widetilde X_{A,C}\to\overline X_{A,C}$ be its minimal desingularization, with lifted zero section $\widetilde\sigma_{A,C}$.  This model need not be relatively-minimal.
Then
\[
\operatorname{Stab}_{G_d^{(3)}(\Fb_q)}(A,C)
\simeq
\Aut_{\Pb^1}
  (\widetilde X_{A,C},\widetilde\sigma_{A,C})
\simeq
\Aut_K(E_{A,C},O).
\]

Moreover, $G_d^{(3)}$ is smooth and connected, and $H^1(\Fb_q,G_d^{(3)})=\{1\}$.  Consequently,
\[
[\cU_d^{(3),0}/G_d^{(3)}](\Fb_q)
\simeq
\cU_d^{(3),0}(\Fb_q)\mathbin{//}G_d^{(3)}(\Fb_q).
\]
\end{prop}

\begin{proof}
De Jong's normalization removes the $a_1$- and $a_3$-terms and, when $j=0$, gives \eqref{eq:char3-normal} \cite[Section~4.13(b2)]{dJ2}.  A global coordinate change between two normalized equations has the form
\[
x\longmapsto u^2x+s,
\qquad
y\longmapsto u^3y+vx+w,
\]
where $u\in\Fb_q^\times$, $s\in V_2(d)$, $v\in V_1(d)$, and $w\in V_3(d)$.  Comparing the $xy$- and $y$-coefficients gives $v=w=0$ because $2\neq0$.  The remaining changes are precisely the elements of $G_d^{(3)}(\Fb_q)$ acting by \eqref{eq:char3-action}.  Thus normalization identifies $\cW_d^{(3),0}(q)$ with the action groupoid $\cU_d^{(3),0}(\Fb_q)\mathbin{//}G_d^{(3)}(\Fb_q)$.

As a scheme, $G_d^{(3)}=\Gb_m\times V_2(d)$, so it is smooth and connected.  Lang's theorem gives $H^1(\Fb_q,G_d^{(3)})=\{1\}$, and Proposition~\ref{prop:rational-points-of-quotients} identifies the action groupoid with $[\cU_d^{(3),0}/G_d^{(3)}](\Fb_q)$.

Let $(X^{\mathrm{rel}},\sigma)\to\Pb^1$ be the relatively-minimal proper regular model of the generic fiber, equipped with its zero section. Restriction to the generic fiber gives
\[
\Aut_{\Pb^1}(X^{\mathrm{rel}},\sigma)
\simeq
\Aut_K(E_{A,C},O).
\]
De Jong's construction obtains $(\widetilde X_{A,C},\widetilde\sigma_{A,C})$ from $(X^{\mathrm{rel}},\sigma)$ by successively blowing up closed points on the successive strict transforms of the zero section \cite[Section~4.11]{dJ2}.  Every automorphism fixes the zero section pointwise and hence preserves each blowup center scheme-theoretically. The universal property of the blowup therefore gives a unique lift through every stage.

Conversely, an automorphism of $(\widetilde X_{A,C},\widetilde\sigma_{A,C})$ restricts to the generic fiber.  The lift of this restriction agrees with the original automorphism on the generic fiber and hence everywhere, since $\widetilde X_{A,C}$ is integral and the target is separated.  The canonical Weierstrass contraction is functorial under isomorphisms of resolved pairs \cite[Sections~4.5--4.8]{dJ2}, so the automorphism also descends uniquely to the Weierstrass surface.  De Jong's coordinate-choice construction identifies the descended automorphism with a global Weierstrass coordinate change preserving the fixed equation.  The coefficient comparison above gives $v=w=0$, so these automorphisms lie in $G_d^{(3)}(\Fb_q)$ and form precisely $\operatorname{Stab}_{G_d^{(3)}(\Fb_q)}(A,C)$.
\end{proof}

For comparison, the weighted cardinality of the entire generically smooth $j=0$ locus of coefficient degree $d$ is
\begin{equation}\label{eq:char3-fulljzero}
\frac{(q^{4d+1}-1)q^{6d+1}}
     {(q-1)q^{2d+1}}
=
q^{4d}\frac{q^{4d+1}-1}{q-1},
\end{equation}
which has leading growth $q^{8d}$.  By contrast, the marked extra-automorphism count computed below has leading growth $q^{6d}$.

\subsection{Classification of fixed pairs}

\begin{lem}
\label{lem:char3-autclassification}
Every origin-preserving $K$-automorphism of the generic fiber $E_{A,C}$ in \eqref{eq:char3-normal} has the form \eqref{eq:char3-change}, with $u\in\Fb_q^\times$ and $s\in V_2(d)$, and satisfies
\begin{equation}\label{eq:char3-autcondition}
u^4=1,
\qquad
s^3+As=(u^6-1)C.
\end{equation}
Apart from the identity $(1,0)$ and inversion $(-1,0)$, the solutions are of the following two types.

\begin{enumerate}[\normalfont(i)]
\item If $u=\pm1$ and $s\neq0$, then $A=-s^2$, while $C$ is arbitrary. The automorphism $(1,s)$ has order three, whereas $(-1,s)$ has order six.

\item If $u^2=-1$, then $A\in V_4(d)\setminus\{0\}$ and $s\in V_2(d)$ are arbitrary, while $C=s^3+As$ is forced.  The automorphism has order four.  The two roots of $u^2=-1$ belong to $\Fb_q$ exactly when $r$ is even.
\end{enumerate}
\end{lem}

\begin{proof}
Proposition~\ref{prop:char3-normalized-quotient} shows that every $K$-defined automorphism of \eqref{eq:char3-normal} is given by \eqref{eq:char3-change}, with $u\in\Fb_q^\times$ and $s\in V_2(d)$. Substitution into the equation gives $u^4A=A$ and $s^3+u^4As=(u^6-1)C$.  Since $A\neq0$, the first equality gives $u^4=1$, and the second then becomes the second equation in \eqref{eq:char3-autcondition}.  Factoring $u^4-1=(u^2-1)(u^2+1)$ shows that either $u=\pm1$ or $u^2=-1$.

If $u=\pm1$, then $u^6=1$, and the second equation in \eqref{eq:char3-autcondition} becomes $s(s^2+A)=0$.  When $s=0$, the choices $u=1$ and $u=-1$ give the identity and inversion, respectively. If $s\neq0$, multiplication by $s$ is injective on global sections because $\Pb^1$ is integral, and hence $A=-s^2$.

If $u^2=-1$, then $u^6=-1$ and $u^6-1=-2=1$ in characteristic three. Thus $C=s^3+As$.

Conversely, every tuple listed in cases~\textup{(i)} and~\textup{(ii)} satisfies \eqref{eq:char3-autcondition}, and hence defines a fixed pair.

Finally, the group law gives
\[
(1,s)^3=(1,0),
\qquad
(-1,s)^3=(-1,0),
\qquad
(u,s)^2=(-1,0)\quad\text{if }u^2=-1.
\]
For $s\neq0$, the first two identities give orders three and six, respectively; the third gives order four.  The element $-1$ is a square in $\Fb_{3^r}$ exactly when $r$ is even.
\end{proof}

\subsection{The ambient marked count}

We now form the corresponding fixed-pair incidence in coefficient degree $d$, before removing the minimality defect and recovering exact Faltings height.

Let
\[
Z_d^{(3)}
=
\left\{
(z,g)\in\cU_d^{(3),0}\times G_d^{(3)}:
g\cdot z=z
\right\}
\]
be the full fixed-pair incidence scheme.  The identity and inversion define the disjoint sections
\[
\mathbf e(A,C)=\bigl((A,C),(1,0)\bigr),
\qquad
\boldsymbol\iota(A,C)=\bigl((A,C),(-1,0)\bigr)
\]
of $Z_d^{(3)}\to\cU_d^{(3),0}$.  Define
\[
Z_d^{(3),\mathrm{ex}}
=
Z_d^{(3)}\setminus
\left(
\mathbf e(\cU_d^{(3),0})
\sqcup
\boldsymbol\iota(\cU_d^{(3),0})
\right).
\]

The group $G_d^{(3)}$ acts on the incidence by change of coordinates and conjugation:
\[
h\cdot(z,g)=\bigl(h\cdot z,hgh^{-1}\bigr).
\]
The elements $(1,0)$ and $(-1,0)$ are central in $G_d^{(3)}$, so the two removed sections are stable and the action restricts to $Z_d^{(3),\mathrm{ex}}$.  Put
\[
Y_d^{(3)}(q)
=
\#_q[Z_d^{(3),\mathrm{ex}}/G_d^{(3)}](\Fb_q).
\]
This is the weighted count of coefficient-degree-$d$ equations marked by an automorphism other than $1$ and $[-1]$, before relative minimality is imposed.

\begin{prop}
\label{prop:char3-ambient-mass}
For every $d\geq0$,
\begin{equation}\label{eq:char3-ambient-mass}
Y_d^{(3)}(q)
=
2P_d^{(3)}(q)+2\epsilon_rQ_d^{(3)}(q),
\end{equation}
where $\epsilon_r$ is the parity indicator defined in
\eqref{eq:parity-indicator}, and
\begin{equation}\label{eq:char3-PQ}
P_d^{(3)}(q)
=
q^{4d}\frac{q^{2d+1}-1}{q-1},
\qquad
Q_d^{(3)}(q)
=
\frac{q^{4d+1}-1}{q-1}.
\end{equation}
\end{prop}

\begin{proof}
Proposition~\ref{prop:char3-normalized-quotient} gives $H^1(\Fb_q,G_d^{(3)})=\{1\}$.  Applying Proposition~\ref{prop:rational-points-of-quotients} to the left $G_d^{(3)}$-scheme $Z_d^{(3),\mathrm{ex}}$ gives
\[
[Z_d^{(3),\mathrm{ex}}/G_d^{(3)}](\Fb_q)
\simeq
Z_d^{(3),\mathrm{ex}}(\Fb_q)\mathbin{//}G_d^{(3)}(\Fb_q).
\]
Orbit--stabilizer therefore gives
\[
Y_d^{(3)}(q)
=
\frac{|Z_d^{(3),\mathrm{ex}}(\Fb_q)|}
{|G_d^{(3)}(\Fb_q)|}.
\]

Write $g=(u,s)\in G_d^{(3)}(\Fb_q)$.  For $h=(v,s_0)\in G_d^{(3)}(\Fb_q)$, \eqref{eq:char3-group-law} gives
\[
hgh^{-1}
=
\bigl(u,v^2s+(1-u^2)s_0\bigr).
\]
Thus conjugation preserves the $u$-coordinate, so each fixed-pair locus with prescribed $u$ is stable and may be counted separately.  Fix $u\in\{1,-1\}$.  By Lemma~\ref{lem:char3-autclassification}, a fixed pair with this value of $u$ is determined by $0\neq s\in V_2(d)$ and $C\in V_6(d)$, with $A=-s^2$.  There are $(q^{2d+1}-1)q^{6d+1}$ such pairs.  Dividing by the group order in \eqref{eq:char3-group-size} gives $P_d^{(3)}(q)$ for each of the two values of $u$.

Now suppose that $r$ is even and fix a root $u\in\Fb_q$ of $u^2=-1$. An order-four fixed pair with this value of $u$ is determined by $0\neq A\in V_4(d)$ and $s\in V_2(d)$, with $C=s^3+As$.  There are $(q^{4d+1}-1)q^{2d+1}$ such pairs.  Division by \eqref{eq:char3-group-size} gives $Q_d^{(3)}(q)$ for each of the two roots.  If $r$ is odd, the equation $u^2=-1$ has no solution in $\Fb_q$, so these loci contribute zero.

The stabilizer of an $\Fb_q$-valued fixed pair $(z,g)$ is $C_{\operatorname{Stab}_{G_d^{(3)}(\Fb_q)}(z)}(g)$. By Proposition~\ref{prop:char3-normalized-quotient}, this is the centralizer of the marked automorphism in the automorphism group of the corresponding resolved model.  Hence orbit--stabilizer gives the required marked-inertia weights.  Adding the two contributions indexed by $u=\pm1$ and, when $r$ is even, the two contributions indexed by $u^2=-1$ proves \eqref{eq:char3-ambient-mass}.
\end{proof}

\section{Extra automorphisms in characteristic 2}
\label{sec:char2-fixed}

Assume throughout that $q=2^r$ with $r\geq1$, and retain the notation $V_i(d)=H^0(\Pb^1,\mathcal O(id))$ for $i,d\geq0$.

\subsection{Normalized equations and the residual coordinate group}

In characteristic two, every generically smooth $j=0$ Weierstrass equation of coefficient degree $d$ admits a normalized form
\begin{equation}\label{eq:char2-normal-form}
E_{A,B,C}\colon y^2+Ay=x^3+Bx+C,
\end{equation}
where $0\neq A\in V_3(d)$, $B\in V_4(d)$, and $C\in V_6(d)$.  Its discriminant is $A^4$, so $A\neq0$ is precisely the condition of generic smoothness.  Put $\cU_d^{(2),0}=(V_3(d)\setminus\{0\})\times V_4(d)\times V_6(d)$.

The \emph{residual coordinate group} consists of the global coordinate changes preserving \eqref{eq:char2-normal-form}.  Its underlying scheme is $G_d^{(2)}=\Gb_m\times V_1(d)\times V_3(d)$, but its multiplication is not the direct-product law.  For $g=(u,s,t)$, define
\begin{equation}\label{eq:char2-coordinate-change}
\Phi_g(x',y')
=
\bigl(u^2x'+s^2,\,u^3y'+u^2sx'+t\bigr).
\end{equation}
The identity $\Phi_{g_1g_2}=\Phi_{g_1}\circ\Phi_{g_2}$ determines the multiplication
\begin{equation}\label{eq:char2-group-law}
\begin{aligned}
&(u_1,s_1,t_1)(u_2,s_2,t_2)\\
&\qquad=
\bigl(u_1u_2,\,u_1s_2+s_1,\,u_1^3t_2+u_1^2s_1s_2^2+t_1\bigr),
\end{aligned}
\end{equation}
with inverse
\begin{equation}\label{eq:char2-group-inverse}
(u,s,t)^{-1}
=
\bigl(u^{-1},\,u^{-1}s,\,u^{-3}(t+s^3)\bigr).
\end{equation}
Pullback by \eqref{eq:char2-coordinate-change} gives the right action
\begin{equation}\label{eq:char2-coefficient-action}
\begin{aligned}
(A,B,C)\mathbin{\triangleleft}(u,s,t)
=
\bigl(&u^{-3}A,\,u^{-4}(B+As+s^4),\\
&u^{-6}(C+Bs^2+s^6+At+t^2)\bigr).
\end{aligned}
\end{equation}
Thus $((A,B,C)\mathbin{\triangleleft}g_1)\mathbin{\triangleleft}g_2=(A,B,C)\mathbin{\triangleleft}(g_1g_2)$.

We write $[\cU_d^{(2),0}/G_d^{(2)}]$ for the quotient associated with the left action $g\cdot z=z\mathbin{\triangleleft}g^{-1}$.  For the fixed-pair calculations below, we retain the equivalent right-action notation because it gives the conditions directly from \eqref{eq:char2-coefficient-action}.  If $z\mathbin{\triangleleft}g=z$, then $g^{-1}\cdot z=z$, so the corresponding marked element for the left action is $g^{-1}$.  An element and its inverse have the same order and the same centralizer, since $C_{G_d^{(2)}}(g)=C_{G_d^{(2)}}(g^{-1})$.

\begin{prop}
\label{prop:char2-normalized-quotient}
Let $\cW_d^{(2),0}(q)$ be the groupoid of generically smooth $j=0$ Weierstrass equations of coefficient degree $d$ over the fixed base $\Pb^1_{\Fb_q}$.  Its morphisms are global Weierstrass coordinate changes over the identity of $\Pb^1_{\Fb_q}$.  Normalization induces an equivalence
\[
\cW_d^{(2),0}(q)
\simeq
[\cU_d^{(2),0}/G_d^{(2)}](\Fb_q).
\]

Let $\overline X_{A,B,C}\to\Pb^1$ be the corresponding coefficient-degree-$d$ Weierstrass surface, and let $\widetilde X_{A,B,C}\to\overline X_{A,B,C}$ be its minimal desingularization, with lifted zero section $\widetilde\sigma_{A,B,C}$.  This model need not be relatively-minimal.  Then
\[
\operatorname{Stab}_{G_d^{(2)}(\Fb_q)}(A,B,C)
\simeq
\Aut_{\Pb^1}\bigl(\widetilde X_{A,B,C},\widetilde\sigma_{A,B,C}\bigr)
\simeq
\Aut_K(E_{A,B,C},O).
\]

Moreover, $G_d^{(2)}$ is smooth and connected, and
\begin{equation}\label{eq:char2-group-size}
|G_d^{(2)}(\Fb_q)|=(q-1)q^{4d+2},
\qquad
H^1(\Fb_q,G_d^{(2)})=\{1\}.
\end{equation}
Consequently,
\[
[\cU_d^{(2),0}/G_d^{(2)}](\Fb_q)
\simeq
\cU_d^{(2),0}(\Fb_q)\mathbin{//}G_d^{(2)}(\Fb_q).
\]
\end{prop}

\begin{proof}
For a generalized Weierstrass equation in characteristic two, the conditions $j=0$ and generic smoothness force $a_1=0$.  De Jong's normalization eliminates the $x^2$ coefficient \cite[Section~4.13(c2)]{dJ2}, giving \eqref{eq:char2-normal-form}.

Every isomorphism of generalized Weierstrass equations is induced by an admissible coordinate change \cite[Sections~4.8--4.9]{dJ2}; see also \cite[Chapter~III, Section~1]{Silverman}.  Between two normalized equations it has the form
\[
x=u^2x'+\rho,
\qquad
y=u^3y'+u^2\eta x'+t,
\]
where $u\in\Fb_q^\times$, $\rho\in V_2(d)$, $\eta\in V_1(d)$, and $t\in V_3(d)$.  Comparing the $x'^2$ coefficients gives $\rho=\eta^2$.  Writing $s=\eta$ yields \eqref{eq:char2-coordinate-change}, while direct composition and substitution give \eqref{eq:char2-group-law}--\eqref{eq:char2-coefficient-action}.  Thus normalization identifies $\cW_d^{(2),0}(q)$ with the action groupoid $\cU_d^{(2),0}(\Fb_q)\mathbin{//}G_d^{(2)}(\Fb_q)$ for the associated left action.

As a scheme, $G_d^{(2)}=\Gb_m\times V_1(d)\times V_3(d)$, so it is smooth and connected.  Since $\dim_{\Fb_q}V_i(d)=id+1$, its group of rational points has order $(q-1)q^{(d+1)+(3d+1)}=(q-1)q^{4d+2}$.  Lang's theorem gives $H^1(\Fb_q,G_d^{(2)})=\{1\}$, and Proposition~\ref{prop:rational-points-of-quotients} identifies this action groupoid with $[\cU_d^{(2),0}/G_d^{(2)}](\Fb_q)$.

Let $(X^{\mathrm{rel}},\sigma)\to\Pb^1$ be the relatively-minimal proper regular model of the generic fiber, equipped with its zero section.  Restriction to the generic fiber gives
\[
\Aut_{\Pb^1}(X^{\mathrm{rel}},\sigma)
\simeq
\Aut_K(E_{A,B,C},O).
\]

The lifting-and-contraction argument in the proof of Proposition~\ref{prop:char3-normalized-quotient} applies without change.  The prescribed blowups are functorial under automorphisms fixing the zero section, and the canonical Weierstrass contraction is functorial under isomorphisms of resolved pairs.  Consequently,
\[
\Aut_{\Pb^1}\bigl(\widetilde X_{A,B,C},\widetilde\sigma_{A,B,C}\bigr)
\simeq
\Aut_{\Pb^1}(X^{\mathrm{rel}},\sigma)
\simeq
\Aut_K(E_{A,B,C},O).
\]
De Jong's coordinate-choice construction identifies the first group with the stabilizer of the fixed normalized equation in the full group of global Weierstrass coordinate changes.  The coefficient comparison above gives $\rho=\eta^2$, so the corresponding coordinate changes lie in $G_d^{(2)}(\Fb_q)$ and are precisely the elements of $\operatorname{Stab}_{G_d^{(2)}(\Fb_q)}(A,B,C)$.
\end{proof}

For comparison, the weighted cardinality of the entire generically smooth $j=0$ locus of coefficient degree $d$ is
\begin{equation}\label{eq:char2-full-jzero}
\frac{(q^{3d+1}-1)q^{4d+1}q^{6d+1}}
{(q-1)q^{4d+2}}
=
q^{6d}\frac{q^{3d+1}-1}{q-1}.
\end{equation}
It has leading growth $q^{9d}$.  By contrast, the marked extra-automorphism count computed below has leading growth $q^{5d}$ when $r$ is odd and $q^{6d}$ when $r$ is even.

\subsection{Classification of fixed pairs}

An $\Fb_q$-valued \emph{fixed pair} is a pair $(z,g)$ satisfying $z\mathbin{\triangleleft}g=z$, where $z=(A,B,C)\in\cU_d^{(2),0}(\Fb_q)$ and $g\in G_d^{(2)}(\Fb_q)$.  Here $\triangleleft$ is the right coefficient action in \eqref{eq:char2-coefficient-action}.  Write $g=(u,s,t)$ and put $v=u^2s$.

\begin{lem}
\label{lem:char2-fixed-pair-classification}
The residual coordinate change $g=(u,s,t)$ stabilizes $E_{A,B,C}$ if and only if
\begin{equation}\label{eq:char2-u-condition}
u^3=1
\end{equation}
and
\begin{align}
v^4+Av+(1+u^2)B&=0,
\label{eq:char2-B-condition}\\
t^2+At&=v^6+u^2Bv^2.
\label{eq:char2-t-condition}
\end{align}
Under these conditions, $s=uv$, and the induced automorphism is
\begin{equation}\label{eq:char2-automorphism-map}
\phi_{u,v,t}(x,y)
=
\bigl(u^2x+u^2v^2,\,y+vx+t\bigr).
\end{equation}
In these coordinates, inversion is the nontrivial involution $[-1](x,y)=(x,y+A)$, represented by $(u,v,t)=(1,0,A)$.  Apart from the identity $(1,0,0)$ and inversion, the $\Fb_q$-defined solutions are as follows.

\begin{enumerate}[(I)]
\item Suppose that $u=1$ and $v\neq0$.  There is a unique $w\in V_2(d)$ such that
\begin{equation}\label{eq:char2-type-I}
A=v^3,
\qquad
t=vw,
\qquad
B=w^2+v^2w+v^4,
\end{equation}
while $C$ is arbitrary.  The marked automorphism has order four and its square is inversion.

\item Suppose that $u\in\mu_3(\Fb_q)\setminus\{1\}$, where $\mu_3(\Fb_q)=\{a\in\Fb_q^\times:a^3=1\}$.  The two nontrivial cube roots of unity belong to $\Fb_q$ if and only if $r$ is even.  For each such $u$, the sections $A\in V_3(d)\setminus\{0\}$ and $v\in V_1(d)$ are arbitrary, while
\begin{equation}\label{eq:char2-type-II-B}
B=u^2(v^4+Av)
\end{equation}
and
\begin{equation}\label{eq:char2-type-II-t}
t=uv^3+\delta A,
\qquad
\delta\in\{0,1\}.
\end{equation}
For $\delta=0$, the marked automorphism has order three.  For $\delta=1$, it has order six and its cube is inversion.  Again, $C$ is arbitrary.
\end{enumerate}
\end{lem}

\begin{proof}
Equating \eqref{eq:char2-coefficient-action} with $(A,B,C)$ gives
\[
(u^3-1)A=0,
\qquad
(u^4-1)B+As+s^4=0,
\qquad
(u^6-1)C+Bs^2+s^6+At+t^2=0.
\]
Since $A\neq0$, the first equation gives $u^3=1$.  Because $v=u^2s$, one then has $s=uv$.  Multiplying the second equation by $u^2$ and using $u^3=1$ gives \eqref{eq:char2-B-condition}; rewriting the third in the same variables gives \eqref{eq:char2-t-condition}.

Suppose first that $u=1$.  If $v=0$, then $t(t+A)=0$.  Since $\Pb^1$ is integral and $A\neq0$, this gives $t=0$ or $t=A$, corresponding respectively to the identity and inversion.

Now suppose that $v\neq0$.  Equation~\eqref{eq:char2-B-condition} becomes $v(v^3+A)=0$.  Since $\Pb^1$ is integral, multiplication by the nonzero section $v$ is injective on global sections, and hence $A=v^3$.

We claim that $v$ divides $t$.  If $v$ is nowhere vanishing, then $d=0$ and the assertion is immediate.  Otherwise, let $x$ be a zero of $v$ of order $a>0$.  If $\operatorname{ord}_x(t)=b<a$, then $A=v^3$ gives
\[
\operatorname{ord}_x(t^2)=2b<3a+b=\operatorname{ord}_x(At).
\]
Hence no cancellation occurs between these terms, so $\operatorname{ord}_x(t^2+At)=2b<2a$, whereas regularity of $B$ at $x$ gives $\operatorname{ord}_x(v^6+Bv^2)\geq2a$, a contradiction.  Thus $\operatorname{ord}_x(t)\geq\operatorname{ord}_x(v)$ at every zero of $v$.  Therefore $t/v$ has no poles and defines a unique global section $w=t/v\in V_2(d)$.  Substituting $t=vw$ into \eqref{eq:char2-t-condition} and dividing by $v^2$ gives \eqref{eq:char2-type-I}.

Now suppose that $u\neq1$.  Since $u^3=1$, one has $u^2+u+1=0$.  Equation~\eqref{eq:char2-B-condition} gives \eqref{eq:char2-type-II-B}, and after this substitution \eqref{eq:char2-t-condition} becomes
\[
(t+uv^3)(t+uv^3+A)=0.
\]
Since $\Pb^1$ is integral, one factor vanishes identically, giving \eqref{eq:char2-type-II-t}.  The two nontrivial cube roots of unity lie in $\Fb_{2^r}$ exactly when $r$ is even.

Conversely, direct substitution into \eqref{eq:char2-u-condition}--\eqref{eq:char2-t-condition} shows that every tuple listed in cases~\textup{(I)} and~\textup{(II)} is a fixed pair.

Finally, direct composition gives
\[
\phi_{1,v,t}^2(x,y)=(x,y+v^3)=(x,y+A)
\]
in type~(I), so the square is inversion.  In type~(II), one obtains
\[
\phi_{u,v,t}^3(x,y)=(x,y+t+uv^3).
\]
When $\delta=0$, this is the identity, and when $\delta=1$, it is inversion.  Since the $u$-coordinate has order three, the marked automorphism has exact order three and six, respectively.
\end{proof}

\subsection{The ambient marked count}

Let
\[
Z_d^{(2)}
=
\left\{(z,g)\in\cU_d^{(2),0}\times G_d^{(2)}:z\mathbin{\triangleleft}g=z\right\}
\]
be the full fixed-pair incidence scheme.  The identity and inversion define the disjoint sections
\[
\mathbf e(A,B,C)=\bigl((A,B,C),(1,0,0)\bigr),
\qquad
\boldsymbol\iota(A,B,C)=\bigl((A,B,C),(1,0,A)\bigr)
\]
of $Z_d^{(2)}\to\cU_d^{(2),0}$.  Define
\[
Z_d^{(2),\mathrm{ex}}
=
Z_d^{(2)}\setminus
\left(\mathbf e(\cU_d^{(2),0})\sqcup\boldsymbol\iota(\cU_d^{(2),0})\right).
\]

The group $G_d^{(2)}$ acts on the incidence by change of coordinates and conjugation:
\[
(z,g)\mathbin{\triangleleft}h
=
\bigl(z\mathbin{\triangleleft}h,\,h^{-1}gh\bigr).
\]
As above, the quotient below is taken for the associated left action $h\cdot(z,g)=(z,g)\mathbin{\triangleleft}h^{-1}=(z\mathbin{\triangleleft}h^{-1},hgh^{-1})$.  The identity section is stable.  If $h=(u,s,t)$, then the group law gives $h^{-1}(1,0,A)h=(1,0,u^{-3}A)$.  Since $u^{-3}A$ is the first coefficient of $z\mathbin{\triangleleft}h$, the inversion section is also stable.  Hence the action restricts to $Z_d^{(2),\mathrm{ex}}$.  Put
\[
Y_d^{(2)}(q)
=
\#_q[Z_d^{(2),\mathrm{ex}}/G_d^{(2)}](\Fb_q).
\]
This is the weighted count of coefficient-degree-$d$ equations marked by an automorphism other than $1$ and $[-1]$, before relative minimality is imposed.

\begin{prop}
\label{prop:char2-ambient-mass}
For every $d\geq0$,
\begin{equation}\label{eq:char2-ambient-mass}
Y_d^{(2)}(q)
=
P_d^{(2)}(q)+4\epsilon_rQ_d^{(2)}(q),
\end{equation}
where $\epsilon_r$ is the parity indicator defined in \eqref{eq:parity-indicator}, and
\begin{equation}\label{eq:char2-PQ}
P_d^{(2)}(q)
=
q^{4d}\frac{q^{d+1}-1}{q-1},
\qquad
Q_d^{(2)}(q)
=
q^{3d}\frac{q^{3d+1}-1}{q-1}.
\end{equation}
\end{prop}

\begin{proof}
By Proposition~\ref{prop:char2-normalized-quotient}, $H^1(\Fb_q,G_d^{(2)})=\{1\}$.  Applying Proposition~\ref{prop:rational-points-of-quotients} to the associated left $G_d^{(2)}$-scheme $Z_d^{(2),\mathrm{ex}}$ gives
\[
[Z_d^{(2),\mathrm{ex}}/G_d^{(2)}](\Fb_q)
\simeq
Z_d^{(2),\mathrm{ex}}(\Fb_q)\mathbin{//}G_d^{(2)}(\Fb_q).
\]
Orbit--stabilizer and \eqref{eq:char2-group-size} therefore give
\[
Y_d^{(2)}(q)
=
\frac{|Z_d^{(2),\mathrm{ex}}(\Fb_q)|}{|G_d^{(2)}(\Fb_q)|}.
\]

By \eqref{eq:char2-group-law}, the map $G_d^{(2)}\to\Gb_m$, $(u,s,t)\mapsto u$, is a group homomorphism.  Since $\Gb_m$ is abelian, conjugation preserves the $u$-coordinate.  In type~\textup{(II)}, Lemma~\ref{lem:char2-fixed-pair-classification} gives orders three and six for $\delta=0$ and $\delta=1$, respectively.  Since conjugation preserves order, it also preserves $\delta$.  Consequently, the type~\textup{(I)} locus and each type~\textup{(II)} locus indexed by $(u,\delta)$ are stable and may be counted separately.

A type~(I) fixed pair is determined by $0\neq v\in V_1(d)$, $w\in V_2(d)$, and $C\in V_6(d)$.  Thus there are $(q^{d+1}-1)q^{2d+1}q^{6d+1}$ such pairs.  Dividing by the group order in \eqref{eq:char2-group-size} gives $P_d^{(2)}(q)$.

Suppose that $r$ is even.  For each $u\in\mu_3(\Fb_q)\setminus\{1\}$ and $\delta\in\{0,1\}$, a type~(II) fixed pair is determined by $0\neq A\in V_3(d)$, $v\in V_1(d)$, and $C\in V_6(d)$; the coefficients $B$ and $t$ are forced by \eqref{eq:char2-type-II-B} and \eqref{eq:char2-type-II-t}.  There are $(q^{3d+1}-1)q^{d+1}q^{6d+1}$ such pairs for each fixed $(u,\delta)$, and division by \eqref{eq:char2-group-size} gives $Q_d^{(2)}(q)$.  There are two choices of $u$ and two of $\delta$, giving four fixed-pair loci and the factor $4\epsilon_r$.  When $r$ is odd, $\mu_3(\Fb_q)=\{1\}$, so no $\Fb_q$-valued type~(II) fixed pairs occur.

For the right action, the stabilizer of an $\Fb_q$-valued fixed pair $(z,g)$ is $C_{\operatorname{Stab}_{G_d^{(2)}(\Fb_q)}(z)}(g)$.  By Proposition~\ref{prop:char2-normalized-quotient}, this is the centralizer of the marked automorphism in the automorphism group of the corresponding resolved model.  Thus the quotient assigns each marked pair the reciprocal of its centralizer order, which is precisely its groupoid weight.  Adding the type~(I) contribution and, when $r$ is even, the four type~(II) contributions proves \eqref{eq:char2-ambient-mass}.
\end{proof}

\section{Marked automorphisms and minimality defect}
\label{sec:defect}

For $p\in\{2,3\}$ and $r\geq1$, let $q=p^r$.  We use a superscript $(p)$ on the normalized fixed-pair spaces, residual coordinate groups, and ambient marked counts constructed in the preceding two sections.

The construction below is the marked, characteristic-two-and-three analogue of the minimality-defect stratification of weighted linear series developed by Bejleri--Park--Satriano.  In their setting, if $\lambda_0,\ldots,\lambda_N$ are the weights and $\kappa=\operatorname{lcm}(\lambda_0,\ldots,\lambda_N)$, then the minimality defect is the size of the quotient obtained by dividing the normalized base profile by $\kappa$ \cite[Definition~4.31]{BPS}.  The normalized base-locus divisor has a canonical decomposition
\[
D_{\mathrm{base}}=\kappa D_{\mathrm{def}}+D_{\mathrm{rem}},
\qquad
\deg D_{\mathrm{def}}=e
\]
\cite[Lemma~4.32]{BPS}.  Over $\Pb^1$, after a locally closed stratification, the defect-$e$ stratum consists of a minimal object of height smaller by $e$ together with an effective divisor in $\Sym^e(\Pb^1)$ \cite[Proposition~6.1 and Corollary~6.2]{BPS}; the corresponding motivic ordinary and inertia convolutions appear in \cite[Lemma~8.6, equations~(6)--(7)]{BPS}.

We follow the same sequence for normalized $j=0$ models in characteristics $2$ and $3$: recover the defect divisor, establish the marked defect decomposition, and invert the resulting divisor convolution.  Because residual translations prevent the representative-level coefficient factorization of \cite[Section~6.1]{BPS} from being imposed invariantly on our normalized coefficients, we formulate the defect intrinsically on the resolved elliptic surface and verify that lifting and contraction preserve the centralizer of the marked automorphism.

A \emph{minimal desingularization} resolves the singularities of a chosen Weierstrass surface, but need not be relatively-minimal as an elliptic fibration.  A curve on an elliptic surface is \emph{vertical} if it is contained in a fiber, and a $(-1)$-curve is a smooth rational curve of self-intersection $-1$.  Relative minimalization contracts the vertical $(-1)$-curves arising from the defect.  We will show that a resolved coefficient-degree-$d$ model is determined by its relatively-minimal model of exact Faltings height $n$ together with an effective divisor $D$ of degree $d-n$.

Let $\Ac_d^{(p),\mathrm{ex}}(q)$ be the essentially finite groupoid of resolved coefficient-degree-$d$ models with a marked extra automorphism.  Its objects are triples $(\widetilde X,\widetilde\sigma,\alpha)$ such that $(\widetilde X,\widetilde\sigma)$ is isomorphic over $\Pb^1_{\Fb_q}$ to the minimal desingularization of a generically smooth normalized $j=0$ Weierstrass equation of coefficient degree $d$, and $\alpha\in\Aut_{\Pb^1}(\widetilde X,\widetilde\sigma)\setminus\{1,[-1]\}$.  Here $[-1]$ is the unique extension of inversion on the generic fiber to the resolved pair.  An arrow $\varphi\colon(\widetilde X,\widetilde\sigma,\alpha)\to(\widetilde X',\widetilde\sigma',\alpha')$ is an isomorphism over the identity of $\Pb^1_{\Fb_q}$ satisfying $\varphi\circ\widetilde\sigma=\widetilde\sigma'$ and $\varphi\alpha\varphi^{-1}=\alpha'$.  The conormal bundle $\widetilde\sigma^*\mathcal O_{\widetilde X}(-\widetilde\sigma)$ has degree $d$.

Let $\cM_n^{\mathrm{ex}}(q)$ be the essentially finite groupoid of triples $(X,\sigma,\alpha)$, where $(X,\sigma)$ is a relatively-minimal proper regular elliptic surface over $\Pb^1_{\Fb_q}$ of exact Faltings height $n$ and $\alpha\in\Aut_{\Pb^1}(X,\sigma)\setminus\{1,[-1]\}$.  An arrow $\varphi\colon(X,\sigma,\alpha)\to(X',\sigma',\alpha')$ is an isomorphism over the identity of $\Pb^1_{\Fb_q}$ satisfying $\varphi\circ\sigma=\sigma'$ and $\varphi\alpha\varphi^{-1}=\alpha'$.  Thus in either groupoid the automorphism group of a marked object is the centralizer of its marking in the automorphism group of the underlying pointed surface.

Restriction to the generic fiber identifies relatively-minimal elliptic surfaces over the fixed base with elliptic curves over $K=\Fb_q(t)$ and preserves automorphism groups \cite[Section~4.2]{dJ2}.  In particular,
\begin{equation}\label{eq:generic-fibre-aut}
\Aut_{\Pb^1}(X,\sigma)\simeq\Aut_K(E,O),
\end{equation}
and therefore $\#_q\cM_n^{\mathrm{ex}}(q)=X_n(q)$.

\begin{lem}\label{lem:resolution-bridge}
For every $p\in\{2,3\}$ and $d\geq0$, there are equivalences
\begin{equation}\label{eq:resolution-bridge}
[Z_d^{(p),\mathrm{ex}}/G_d^{(p)}](\Fb_q)
\simeq
Z_d^{(p),\mathrm{ex}}(\Fb_q)\mathbin{//}G_d^{(p)}(\Fb_q)
\simeq
\Ac_d^{(p),\mathrm{ex}}(q).
\end{equation}
Consequently, $Y_d^{(p)}(q)=\#_q\Ac_d^{(p),\mathrm{ex}}(q)$.
\end{lem}

\begin{proof}
The explicit descriptions of the residual coordinate groups show that $G_d^{(p)}$ is smooth and connected.  Lang's theorem gives $H^1(\Fb_q,G_d^{(p)})=\{1\}$, so Proposition~\ref{prop:rational-points-of-quotients} gives the first equivalence in \eqref{eq:resolution-bridge}.

Propositions~\ref{prop:char3-normalized-quotient} and \ref{prop:char2-normalized-quotient} identify the underlying action groupoids with the groupoids of normalized coefficient-degree-$d$ equations.  Passing to the fixed-pair incidence and deleting the identity and inversion sections identifies $Z_d^{(p),\mathrm{ex}}(\Fb_q)\mathbin{//}G_d^{(p)}(\Fb_q)$ with the groupoid of such equations marked by an extra automorphism.  Minimal desingularization therefore defines a functor to $\Ac_d^{(p),\mathrm{ex}}(q)$.

A coordinate change between Weierstrass equations lifts uniquely to their minimal desingularizations.  Conversely, the \emph{canonical Weierstrass model} of a resolved pair is obtained by contracting in every fiber the irreducible components disjoint from the zero section; this construction is functorial and recovers the original model when the resolved pair is its minimal desingularization.  Hence an isomorphism of resolved pairs descends to an isomorphism of their canonical Weierstrass models \cite[Sections~4.7--4.9]{dJ2}, and the marked automorphism of any object of $\Ac_d^{(p),\mathrm{ex}}(q)$ descends as well.  After choosing normalized presentations, Propositions~\ref{prop:char3-normalized-quotient} and \ref{prop:char2-normalized-quotient} identify the descended automorphisms and isomorphisms with elements of $G_d^{(p)}(\Fb_q)$.  The functor is therefore essentially surjective and fully faithful.

These constructions preserve the marked automorphism and carry conjugation to conjugation.  They therefore identify the corresponding centralizers and preserve weighted cardinality.
\end{proof}

Recall that
\[
c_e(q)=|\Sym^e(\Pb^1)(\Fb_q)|=\frac{q^{e+1}-1}{q-1}.
\]
An $\Fb_q$-point of $\Sym^e(\Pb^1)$ is equivalently a Frobenius-stable effective divisor $D=\sum_{x\in|\Pb^1|}e_x\,x$ of degree $\deg D=\sum_xe_x\deg(x)=e$.  This is the divisor in de Jong's blowup description: each blowup above a closed point $x$ increases the degree of the section conormal bundle by $\deg(x)$.  By analogy with \cite[Definition~4.31]{BPS}, we call $e=d-n=\deg D$ the \emph{minimality defect}.  More precisely, $D$ plays the role of the effective quotient divisor in \cite[Lemma~4.32]{BPS}; it is intrinsic rather than determined by visible divisibility of the coefficients in a chosen normalized equation.

The map $\psi_{n,e}$ of \cite[Section~6.1 and Proposition~6.1]{BPS} is written on representatives as
\[
(f_0,\ldots,f_N)=(g_0h^{\lambda_0},\ldots,g_Nh^{\lambda_N}),
\]
where $h$ is a degree-$e$ form.  The following lemma is the intrinsic resolved-surface analogue of this factorization.

\begin{lem}
\label{lem:blowup-weierstrass}
Let $f\colon(X,\sigma)\to\Pb^1_{\Fb_q}$ be a relatively-minimal proper regular elliptic surface of exact Faltings height $n$, write $\mathcal L_f=(R^1f_*\mathcal O_X)^\vee$ for its fundamental line bundle, and let $D=\sum_{x\in|\Pb^1|}e_x\,x$ be an effective divisor.  Put $d=n+\deg D$.  For each closed point $x$, blow up $\sigma(x)$ and then the closed point above $x$ on each successive strict transform of the section, a total of $e_x$ times, and denote the resulting pair by $(\widetilde X_D,\widetilde\sigma_D)$.  Then
\begin{equation}\label{eq:blowup-conormal}
\widetilde\sigma_D^*\mathcal O_{\widetilde X_D}(-\widetilde\sigma_D)
\simeq
\mathcal L_f\otimes\mathcal O_{\Pb^1}(D),
\end{equation}
so this line bundle has degree $d$.  The canonical Weierstrass model associated to $(\widetilde X_D,\widetilde\sigma_D)$ has coefficient degree $d$, and $\widetilde X_D$ is its minimal desingularization.

The construction is functorial under isomorphisms over the identity of $\Pb^1_{\Fb_q}$.  In particular, lifting induces an isomorphism
\[
\Aut_{\Pb^1}(X,\sigma)
\simeq
\Aut_{\Pb^1}(\widetilde X_D,\widetilde\sigma_D)
\]
that carries $1$ and $[-1]$ to their corresponding lifts, respects conjugation, and identifies the centralizers of corresponding marked automorphisms.

If $p\in\{2,3\}$ and the generic fiber has $j=0$, the corresponding normalization preserves the isomorphism class of the resolved pair and the centralizer of a marked automorphism.
\end{lem}

\begin{proof}
The section self-intersection formula gives $\sigma^*\mathcal O_X(-\sigma)\simeq\mathcal L_f$, with $\deg\mathcal L_f=n$.  Blowing up the section above a closed point $x$ tensors the conormal bundle of its strict transform by $\mathcal O_{\Pb^1}(x)$.  Iterating over $D$ proves \eqref{eq:blowup-conormal}, whose degree is $n+\deg D=d$.

We verify the minimal-resolution criterion of \cite[Section~4.7]{dJ2}.  After base change to $\overline{\Fb}_q$, fix a geometric point $\bar x$ above a closed point $x$ with $e_x>0$.  The zero section lies in the smooth locus of $f$, so $\sigma(\bar x)$ is a smooth point of the fiber.  Let $E_{\bar x,1},\ldots,E_{\bar x,e_x}$ be the successive exceptional curves.  For $1\leq i<e_x$, the curve $E_{\bar x,i}$ is blown up exactly once more at its intersection with the successive strict transform of the section, so its final strict transform is a vertical $(-2)$-curve disjoint from the section.  The terminal curve $E_{\bar x,e_x}$ is a vertical $(-1)$-curve meeting $\widetilde\sigma_D$.

No strict transform of a vertical curve on $X_{\overline{\Fb}_q}$ produces another $(-1)$-curve.  Every smooth rational component of a reducible fiber has self-intersection at most $-2$ by relative minimality, and blowing up can only decrease it.  If the fiber through $\bar x$ is irreducible, its strict transform after the first blowup has self-intersection $-1$, but retains arithmetic genus one and hence is not a $(-1)$-curve.  All later centers above $\bar x$ lie on exceptional curves.  Thus the terminal exceptional curves are precisely the vertical $(-1)$-curves of $\widetilde X_D$, and all meet the section.

The criterion of \cite[Section~4.7]{dJ2} now shows that $\widetilde X_D$ is the minimal desingularization of its canonical Weierstrass contraction.  De Jong's coordinate construction \cite[Sections~4.5--4.6]{dJ2} writes this model with coefficients
\[
a_i\in H^0\bigl(\Pb^1,(\widetilde\sigma_D^*\mathcal O_{\widetilde X_D}(-\widetilde\sigma_D))^{\otimes i}\bigr).
\]
Every degree-$d$ line bundle on $\Pb^1_{\Fb_q}$ is isomorphic to $\mathcal O(d)$, so the model has coefficient degree $d$.

Let $\varphi\colon(X,\sigma)\to(X',\sigma')$ be an isomorphism over the identity of $\Pb^1_{\Fb_q}$, and use the same divisor $D$ on both sides.  Since $\varphi\circ\sigma=\sigma'$, it carries every closed and infinitely near blowup center scheme-theoretically to the corresponding center.  The universal property of the blowup gives a unique lift $\widetilde\varphi_D\colon(\widetilde X_D,\widetilde\sigma_D)\to(\widetilde X'_D,\widetilde\sigma'_D)$, compatibly with identities and composition.

Conversely, after base change to $\overline{\Fb}_q$, the terminal exceptional curves are characterized intrinsically as the vertical $(-1)$-curves of the blown-up surface.  They form Frobenius-stable disjoint collections indexed by the closed points $x$ with $e_x>0$.  The simultaneous contraction of each collection is unique, hence Frobenius equivariant, and therefore descends to $\Fb_q$.  Any isomorphism of the blown-up pairs preserves these collections and descends uniquely through their contractions.  After the terminal collection above each $x$ is contracted, the preceding exceptional curves for those $x$ with $e_x\geq2$ become the new terminal geometric $(-1)$-curves.  Iteration recovers $(X,\sigma)$ and all multiplicities $e_x$, so lifting and contraction are mutually inverse on isomorphisms.

The resulting isomorphism of automorphism groups respects composition and conjugation and carries $1$ and $[-1]$ to the corresponding elements; it therefore identifies the stated centralizers.

Finally, suppose $p\in\{2,3\}$ and the generic fiber has $j=0$.  The normalizations in Propositions~\ref{prop:char3-normalized-quotient} and \ref{prop:char2-normalized-quotient} are effected by global admissible Weierstrass coordinate changes.  Such a change induces an isomorphism of the associated Weierstrass models and hence, by uniqueness of minimal desingularizations, an isomorphism of their resolved pairs.  Transport through this isomorphism carries a marked automorphism to the corresponding marked automorphism and preserves its centralizer.
\end{proof}

The preceding lemma yields the marked defect decomposition.

\begin{prop}\label{prop:marked-defect}
For every $p\in\{2,3\}$ and $d\geq0$, relative minimalization together with recovery of the defect divisor induces an equivalence of essentially finite groupoids
\begin{equation}\label{eq:defect-equivalence}
\Ac_d^{(p),\mathrm{ex}}(q)
\simeq
\bigsqcup_{n=0}^{d}
\bigsqcup_{D\in\Sym^{d-n}(\Pb^1)(\Fb_q)}
\cM_n^{\mathrm{ex}}(q),
\end{equation}
where each symmetric-power set is regarded as a discrete groupoid, meaning that it has only identity arrows.  This equivalence identifies the automorphism groups of the underlying pointed surfaces and hence the centralizers of corresponding marked automorphisms.  Consequently,
\begin{equation}\label{eq:defect-recursion}
Y_d^{(p)}(q)
=
\sum_{n=0}^{d}c_{d-n}(q)X_n(q)
=
\sum_{e=0}^{d}c_e(q)X_{d-e}(q).
\end{equation}
\end{prop}

\begin{proof}
Let $(\widetilde X,\widetilde\sigma,\widetilde\alpha)\in\Ac_d^{(p),\mathrm{ex}}(q)$.  By de Jong's description of nonminimal Weierstrass models, there are a relatively-minimal proper regular elliptic surface $(X,\sigma)$ of exact Faltings height $n\leq d$ and an effective divisor $D=\sum_{x\in|\Pb^1|}e_x\,x$ of degree $d-n$ such that $(\widetilde X,\widetilde\sigma)$ is isomorphic to the prescribed blowup of $(X,\sigma)$, with $e_x$ successive blowups above $x$ \cite[Section~4.11]{dJ2}.  Lemma~\ref{lem:blowup-weierstrass} shows that successively contracting the vertical $(-1)$-curves recovers $(X,\sigma)$ and $D$ uniquely.  The marked automorphism $\widetilde\alpha$ therefore descends uniquely to an automorphism $\alpha\in\Aut_{\Pb^1}(X,\sigma)$.  Since the generic fiber is unchanged, $\alpha$ remains distinct from $1$ and $[-1]$.  This construction is functorial on arrows and defines a functor from the left side of \eqref{eq:defect-equivalence} to the right.

Conversely, begin with $(X,\sigma,\alpha)\in\cM_n^{\mathrm{ex}}(q)$ and $D\in\Sym^{d-n}(\Pb^1)(\Fb_q)$.  Lemma~\ref{lem:blowup-weierstrass} constructs the prescribed blowup $(\widetilde X_D,\widetilde\sigma_D)$, the unique lift $\widetilde\alpha_D$, and a coefficient-degree-$d$ Weierstrass model whose minimal desingularization is $(\widetilde X_D,\widetilde\sigma_D)$.  By the classification of automorphism groups of elliptic curves \cite[Appendix~A, Proposition~1.2(c)]{Silverman}, an origin-preserving automorphism other than $1$ and $[-1]$ forces the generic fiber to have $j=0$ in characteristics $2$ and $3$.  The final assertion of Lemma~\ref{lem:blowup-weierstrass} therefore shows that $(\widetilde X_D,\widetilde\sigma_D,\widetilde\alpha_D)$ is an object of $\Ac_d^{(p),\mathrm{ex}}(q)$.

By the preceding construction and Lemma~\ref{lem:blowup-weierstrass}, applying the blowup and contraction procedures successively recovers each object up to canonical isomorphism, and the two procedures induce inverse bijections on morphism sets.  They recover the support and multiplicities of $D$, identify the automorphism groups of the underlying pointed surfaces, and hence identify the centralizers of corresponding marked automorphisms.  An arrow therefore forces the source and target divisors to agree, so there are no arrows between summands indexed by distinct divisors.  This proves \eqref{eq:defect-equivalence}.

Taking weighted cardinalities and using Lemma~\ref{lem:resolution-bridge} and $\#_q\cM_n^{\mathrm{ex}}(q)=X_n(q)$ gives
\[
Y_d^{(p)}(q)
=
\sum_{n=0}^{d}c_{d-n}(q)X_n(q)
=
\sum_{e=0}^{d}c_e(q)X_{d-e}(q),
\]
which is \eqref{eq:defect-recursion}.
\end{proof}

Proposition~\ref{prop:marked-defect} is an equivalence of essentially finite groupoids over $\Fb_q$ and the fixed base $\Pb^1$; it does not assert a family-level equivalence over arbitrary test schemes.  In the weighted-projective setting, the quotient and remainder divisors are defined in families with fixed normalized base profile, and the factorization map becomes an isomorphism after stratifying its source and target \cite[Lemma~4.32, Proposition~6.1, and Corollary~6.2]{BPS}.

The intrinsic formulation is necessary because the representative-level coefficient division used in \cite[Proposition~6.1]{BPS} cannot be imposed invariantly on our normalized representatives.  In characteristic three, a residual translation of the $x$-coordinate can destroy visible divisibility of $C$ without changing the resolved nonminimal model.  Thus, even if $h\in H^0(\Pb^1,\mathcal O(d-n))$ has zero divisor $D$, one cannot require $C=h^6C_0$ in every normalized representative.  In characteristic two, the residual translation parameters $s$ and $t$ mix $B$ with $As+s^4$ and $C$ with $Bs^2+s^6+At+t^2$.  In both characteristics, the blowup divisor is intrinsic, whereas visible coefficient divisibility is not.

\subsection{Inverting the divisor convolution}

The following generating-series argument is the weighted-point-count analogue of the ordinary and inertia factorizations in \cite[Lemma~8.6, equations~(6)--(7)]{BPS}.

\begin{lem}
\label{lem:defect-inversion}
For every $p\in\{2,3\}$ and $M\geq1$,
\begin{equation}
\label{eq:cumulative-defect-inversion}
\sum_{n=0}^{M}X_n(q)
=
Y_M^{(p)}(q)-qY_{M-1}^{(p)}(q).
\end{equation}
\end{lem}

\begin{proof}
Working with formal power series with rational coefficients, put
\[
C(z)=\sum_{e\geq0}c_e(q)z^e=\frac{1}{(1-z)(1-qz)},
\qquad
X(z)=\sum_{n\geq0}X_n(q)z^n,
\qquad
Y^{(p)}(z)=\sum_{d\geq0}Y_d^{(p)}(q)z^d.
\]
Equation~\eqref{eq:defect-recursion} gives $Y^{(p)}(z)=C(z)X(z)$.  Consequently,
\[
\sum_{M\geq0}\left(\sum_{n=0}^{M}X_n(q)\right)z^M
=
\frac{X(z)}{1-z}
=
(1-qz)Y^{(p)}(z).
\]
Taking the coefficient of $z^M$ proves \eqref{eq:cumulative-defect-inversion}.
\end{proof}

\subsection{The characteristic-three specialization}

For $q=3^r$, the ambient fixed-pair calculation gives
\begin{equation}\label{eq:char3-ambient-decomposition}
Y_d^{(3)}(q)=2P_d^{(3)}(q)+2\epsilon_rQ_d^{(3)}(q),
\end{equation}
where $P_d^{(3)}(q)=q^{4d}(q^{2d+1}-1)/(q-1)$ and $Q_d^{(3)}(q)=(q^{4d+1}-1)/(q-1)$.  For every $M\geq1$, direct subtraction using \eqref{eq:alpha-beta} gives
\begin{align*}
P_M^{(3)}(q)-qP_{M-1}^{(3)}(q)
&=
\alpha_6(q)q^{6M}-\beta_4(q)q^{4M},\\
Q_M^{(3)}(q)-qQ_{M-1}^{(3)}(q)
&=
\alpha_4(q)q^{4M}+1.
\end{align*}

\begin{cor}
\label{cor:char3-extra-inertia-sum}
For every $M\geq1$,
\begin{align}
\sum_{n=0}^{M}X_n(q)
={}&
2\bigl(\alpha_6(q)q^{6M}-\beta_4(q)q^{4M}\bigr)
\nonumber\\
&+2\epsilon_r\bigl(\alpha_4(q)q^{4M}+1\bigr).
\label{eq:char3-extra-inertia-sum}
\end{align}
\end{cor}

\begin{proof}
By Lemma~\ref{lem:defect-inversion} and \eqref{eq:char3-ambient-decomposition},
\[
\sum_{n=0}^{M}X_n(q)
=
2\bigl(P_M^{(3)}(q)-qP_{M-1}^{(3)}(q)\bigr)
+
2\epsilon_r\bigl(Q_M^{(3)}(q)-qQ_{M-1}^{(3)}(q)\bigr).
\]
Substitution of the preceding two identities proves the result.
\end{proof}

\subsection{The characteristic-two specialization}

For $q=2^r$, the ambient fixed-pair calculation gives
\begin{equation}\label{eq:char2-ambient-decomposition}
Y_d^{(2)}(q)=P_d^{(2)}(q)+4\epsilon_rQ_d^{(2)}(q),
\end{equation}
where $P_d^{(2)}(q)=q^{4d}(q^{d+1}-1)/(q-1)$ and $Q_d^{(2)}(q)=q^{3d}(q^{3d+1}-1)/(q-1)$.  For every $M\geq1$, direct subtraction using \eqref{eq:alpha-beta} gives
\begin{align*}
P_M^{(2)}(q)-qP_{M-1}^{(2)}(q)
&=
\alpha_5(q)q^{5M}-\beta_4(q)q^{4M},\\
Q_M^{(2)}(q)-qQ_{M-1}^{(2)}(q)
&=
\alpha_6(q)q^{6M}-\beta_3(q)q^{3M}.
\end{align*}

\begin{cor}
\label{cor:char2-extra-inertia-sum}
For every $M\geq1$,
\begin{align}
\sum_{n=0}^{M}X_n(q)
={}&
\alpha_5(q)q^{5M}-\beta_4(q)q^{4M}
\nonumber\\
&+4\epsilon_r\bigl(\alpha_6(q)q^{6M}-\beta_3(q)q^{3M}\bigr).
\label{eq:char2-extra-inertia-sum}
\end{align}
\end{cor}

\begin{proof}
By Lemma~\ref{lem:defect-inversion} and \eqref{eq:char2-ambient-decomposition},
\[
\sum_{n=0}^{M}X_n(q)
=
P_M^{(2)}(q)-qP_{M-1}^{(2)}(q)
+
4\epsilon_r\bigl(Q_M^{(2)}(q)-qQ_{M-1}^{(2)}(q)\bigr).
\]
Substitution of the preceding two identities proves the result.
\end{proof}

At $d=0$ no positive defect is possible, so
\begin{equation}
\label{eq:defect-height-zero}
X_0(q)=Y_0^{(p)}(q)
=
\begin{cases}
2+2\epsilon_r,&q=3^r,\\
1+4\epsilon_r,&q=2^r.
\end{cases}
\end{equation}
These are the omitted $M=0$ values; the characteristic-specific closed formulas above require $M\geq1$.

\section{Proofs of the counting formulas}
\label{sec:count}

\begin{thm}
\label{thm:faltings-height}
Let $p\in\{2,3\}$ and $r\geq1$, and put $q=p^r$ and $K=\Fb_q(t)$.  For every $M\in\Zb_{\geq1}$, the following formulas hold.

If $p=3$, then
\begin{align*}
\Nc(K,q^{12M})
={}&
2\left(
\frac{q^9-1}{q^8-q^7}q^{10M}
-q^{2M}
-\beta_4(q)q^{4M}
-1
\right)\\
&+
2\left(
\alpha_6(q)q^{6M}
-\beta_4(q)q^{4M}
\right)
+
2\epsilon_r\left(
\alpha_4(q)q^{4M}+1
\right).
\end{align*}

If $p=2$, then
\begin{align*}
\Nc(K,q^{12M})
={}&
2\left(
\frac{q^9-1}{q^8-q^7}q^{10M}
-q^{2M}
-\beta_6(q)q^{6M}
\right)\\
&+
\alpha_5(q)q^{5M}
-\beta_4(q)q^{4M}
+
4\epsilon_r\left(
\alpha_6(q)q^{6M}
-\beta_3(q)q^{3M}
\right).
\end{align*}

For $M=0$,
\[
\Nc(K,1)
=
\begin{cases}
2q+2+2\epsilon_r,&p=3,\\
2q+1+4\epsilon_r,&p=2.
\end{cases}
\]
Equivalently,
\[
\Nc(K,1)
=
\begin{cases}
2q+2,&p=3\text{ and }r\text{ odd},\\
2q+4,&p=3\text{ and }r\text{ even},\\
2q+1,&p=2\text{ and }r\text{ odd},\\
2q+5,&p=2\text{ and }r\text{ even}.
\end{cases}
\]
\end{thm}

\begin{proof}
Since the Faltings height is a nonnegative integer and $U_n(q)$ counts curves of exact height $n$, one has
\[
\Nc(K,q^{12M})
=
\sum_{n=0}^M U_n(q)
\]
for every $M\geq0$.  For $M\geq1$, summing the marked-inertia identity $U_n(q)=2I_n(q)+X_n(q)$ from Lemma~\ref{lem:inertia} gives
\[
\Nc(K,q^{12M})
=
2\sum_{n=0}^M I_n(q)
+
\sum_{n=0}^M X_n(q).
\]
Theorem~\ref{thm:weighted} evaluates the first sum, while Corollaries~\ref{cor:char3-extra-inertia-sum} and \ref{cor:char2-extra-inertia-sum} evaluate the second in characteristics three and two, respectively.  Substitution gives the asserted positive-height formulas.

For $M=0$, Theorem~\ref{thm:weighted} gives $I_0(q)=q$, while \eqref{eq:defect-height-zero} gives the corresponding value of $X_0(q)$.  The marked-inertia identity therefore yields
\[
\Nc(K,1)=U_0(q)
=
\begin{cases}
2q+2+2\epsilon_r,&q=3^r,\\
2q+1+4\epsilon_r,&q=2^r.
\end{cases}
\]
The equivalent four-case formulation follows from the definition of $\epsilon_r$.
\end{proof}

\begin{proof}[Proof of Theorem~\ref{thm:main}]
For real $B\geq1$, let $B^\ast$ be the largest element of $q^{12\Zb_{\geq0}}$ not exceeding $B$.  Since $H_\Delta(E)=q^{12\height(E)}$ and $\height(E)\in\Zb_{\geq0}$, one has $H_\Delta(E)\leq B$ if and only if $H_\Delta(E)\leq B^\ast$.  Consequently,
\[
\Nc^w(K,B)=\Nc^w(K,B^\ast),
\qquad
\Nc(K,B)=\Nc(K,B^\ast).
\]

If $B^\ast=1$, then $\Nc^w(K,B)=\Nc^w(K,1)=I_0(q)=q$ by Theorem~\ref{thm:weighted}, while $\Nc(K,B)=\Nc(K,1)$ is given by Theorem~\ref{thm:faltings-height}.  This proves the height-zero formulas.

Suppose that $B^\ast>1$ and write $B^\ast=q^{12M}$ with $M\geq1$.  Then
\[
\begin{gathered}
q^{10M}=(B^\ast)^{5/6},
\qquad
q^{6M}=(B^\ast)^{1/2},
\qquad
q^{5M}=(B^\ast)^{5/12},\\
q^{4M}=(B^\ast)^{1/3},
\qquad
q^{3M}=(B^\ast)^{1/4},
\qquad
q^{2M}=(B^\ast)^{1/6}.
\end{gathered}
\]
By definition,
\[
\Nc^w(K,B^\ast)
=
\sum_{n=0}^M I_n(q),
\qquad
\Nc(K,B^\ast)
=
\sum_{n=0}^M U_n(q).
\]
Theorem~\ref{thm:weighted} evaluates the first sum, and Theorem~\ref{thm:faltings-height} evaluates the second.  Substituting the displayed identities gives the corresponding formulas for $\Nc^w(K,B^\ast)$ and $\Nc(K,B^\ast)$; the initial equalities then complete the proof for arbitrary $B$.
\end{proof}

\section{Twists and the lower-order terms}
\label{sec:twists}

\subsection{Finite-field twists at height zero}

By an \emph{$\Fb_q$-form}, or \emph{twist}, of an elliptic curve $E/\overline{\Fb}_q$ we mean an elliptic curve $E'/\Fb_q$ such that $E'_{\overline{\Fb}_q}\simeq E$.  Kronberg--Soomro--Top determine the finite-field forms of elliptic curves in characteristics $2$ and $3$.  Since all elliptic curves with $j=0$ are isomorphic over $\overline{\Fb}_q$ in either characteristic, their results give
\[
\#\left\{
\begin{array}{c}
\Fb_q\text{-isomorphism classes}\\
\text{of elliptic curves with }j=0
\end{array}
\right\}
=
\begin{cases}
4+2\epsilon_r,&q=3^r,\\
3+4\epsilon_r,&q=2^r
\end{cases},
\]
where the characteristic-three statement is \cite[Proposition~2.2]{KST} and the characteristic-two statement is \cite[Proposition~3.2]{KST}.

Write $I_0^{j=0}(q)$, $X_0^{j=0}(q)$, and $U_0^{j=0}(q)$ for the weighted, extra-inertia, and unweighted height-zero counts restricted to the locus $j=0$.

Every height-zero elliptic curve over $K=\Fb_q(t)$ is constant \cite[Section~4.15]{dJ2}.  Indeed, its fundamental line bundle has degree zero and hence is isomorphic to $\mathcal O_{\Pb^1}$.  After choosing such a trivialization, the coefficients of a global minimal Weierstrass equation lie in $H^0(\Pb^1,\mathcal O_{\Pb^1})=\Fb_q$, so its generic fiber is obtained by base change from $\Fb_q$.  Conversely, every elliptic curve over $\Fb_q$ has height zero after base change to $K$.

Base change is also injective on isomorphism classes.  Indeed, if two elliptic curves over $\Fb_q$ become isomorphic over $K$, their finite \'etale $\Fb_q$-scheme of origin-preserving isomorphisms has a $K$-point.  This point factors through a closed point whose residue field is a finite extension of $\Fb_q$ embedded in $K$.  Since $\Fb_q$ is relatively algebraically closed in $K$, that residue field is $\Fb_q$.  Thus the curves are already isomorphic over $\Fb_q$, and base change induces a bijection between $\Fb_q$-isomorphism classes and height-zero $K$-isomorphism classes.

At coefficient degree $d=0$, no positive minimality defect is possible.  Setting $d=0$ in \eqref{eq:char3-fulljzero} or \eqref{eq:char2-full-jzero} therefore gives $I_0^{j=0}(q)=1$.  By the classification of origin-preserving automorphism groups of elliptic curves \cite[Appendix~A, Proposition~1.2(c)]{Silverman}, automorphisms other than $1$ and $[-1]$ occur only at $j=0$ in characteristics $2$ and $3$.  Hence $X_0^{j=0}(q)=X_0(q)$, and \eqref{eq:defect-height-zero} gives
\[
X_0^{j=0}(q)
=
\begin{cases}
2+2\epsilon_r,&q=3^r,\\
1+4\epsilon_r,&q=2^r.
\end{cases}
\]
Consequently, the marked-inertia identity gives
\[
U_0^{j=0}(q)
=
2I_0^{j=0}(q)+X_0^{j=0}(q)
=
\begin{cases}
4+2\epsilon_r,&q=3^r,\\
3+4\epsilon_r,&q=2^r.
\end{cases}
\]
Thus the height-zero calculation recovers the twist numbers of Kronberg--Soomro--Top.  For each $j\in\Fb_q^\times$, there are exactly two $\Fb_q$-isomorphism classes of elliptic curves with $j$-invariant $j$ \cite[Propositions~2.2 and~3.2]{KST}.  Hence the preceding bijection gives $\Nc(K,1)=2(q-1)+(4+2\epsilon_r)=2q+2+2\epsilon_r$ in characteristic three and $\Nc(K,1)=2(q-1)+(3+4\epsilon_r)=2q+1+4\epsilon_r$ in characteristic two, as asserted in Theorem~\ref{thm:faltings-height}. 

This comparison with the finite-field twist classification is specific to height zero.  Positive-height curves on the $j=0$ locus are nonconstant isotrivial $K$-forms and hence are not classified by $\Fb_q$-forms.  Their $K$-defined extra automorphisms are handled directly by the fixed-pair classifications in Lemmas~\ref{lem:char3-autclassification} and \ref{lem:char2-fixed-pair-classification}, together with the centralizer-weighted counts in Propositions~\ref{prop:char3-ambient-mass} and \ref{prop:char2-ambient-mass}.

\subsection{The characteristic-three lower-order terms}

Suppose that $q=3^r$, and for $M\geq1$ put $B=q^{12M}$.

The dominant marked fixed-pair loci are given by $A=-s^2$ with $s\neq0$, with one locus for each of $u=1$ and $u=-1$.  The sublocus with $u=1$ supports the order-three automorphisms $(x,y)\mapsto(x+s,y)$, while the sublocus with $u=-1$ consists of their products with inversion and supports automorphisms of order six.  Together these subloci have ambient marked count $2P_d^{(3)}(q)$, of leading growth $q^{6d}$.  After inversion of the defect recursion, their combined contribution is
\[
2\bigl(\alpha_6(q)q^{6M}-\beta_4(q)q^{4M}\bigr)
=
2\alpha_6(q)B^{1/2}-2\beta_4(q)B^{1/3}.
\]

The second family of marked fixed-pair loci is cut out by $u^2=-1$.  Their union is defined over $\Fb_q$ for every $r$ and, over $\overline{\Fb}_q$, is the disjoint union of two fixed-pair loci, one for each root of $u^2=-1$.  These loci are individually defined over $\Fb_q$ exactly when $r$ is even.  When $r$ is odd, Frobenius exchanges them, so their union has no $\Fb_q$-points: an $\Fb_q$-point fixed by Frobenius would have to lie in both disjoint geometric components.  The corresponding marked automorphisms have order four, and their contribution is
\[
2\epsilon_r\bigl(\alpha_4(q)q^{4M}+1\bigr)
=
2\epsilon_r\alpha_4(q)B^{1/3}+2\epsilon_r.
\]

The negative terms of order $B^{1/3}$ have two distinct origins.  The summand $-2\beta_4(q)B^{1/3}$ from the marked loci is created when the minimality defect is removed from the positive ambient marked count.  When $r$ is odd, it is the entire extra-inertia contribution of order $B^{1/3}$.  A second contribution $-2\beta_4(q)B^{1/3}-2$ is the correction term in $2\Nc^w(K,B)$ arising from the additional nonsmooth equations in the quasi-elliptic-type stratum of Proposition~\ref{prop:nonsmooth}.  Together these contributions explain the coefficient $2\bigl(\epsilon_r\alpha_4(q)-2\beta_4(q)\bigr)$ and the constant term $2\epsilon_r-2$ in Theorem~\ref{thm:main}.

\subsection{The characteristic-two lower-order terms}

Suppose that $q=2^r$, and for $M\geq1$ put $B=q^{12M}$.  Over an algebraic closure, the origin-preserving automorphism group of a $j=0$ elliptic curve in characteristic two has order $24$.  Its order-four elements form one conjugacy class of six elements, with centralizer of order four.  Its order-three elements form two conjugacy classes of four elements each, and its order-six elements likewise form two conjugacy classes of four elements each; their centralizers have order six \cite[Proposition~3.1]{KST}.

Over $\overline{\Fb}_q$, the geometric type-(II) fixed pairs are partitioned into four loci indexed by $u\in\mu_3(\overline{\Fb}_q)\setminus\{1\}$ and $\delta\in\{0,1\}$.  These loci are individually defined over $\Fb_q$ precisely when $\mu_3(\overline{\Fb}_q)\subset\Fb_q^\times$, equivalently when $\Fb_4\subset\Fb_q$, or equivalently when $r$ is even.  When $r$ is odd, Frobenius exchanges the two loci indexed by the nontrivial cube roots of unity for each value of $\delta$.  Since the exchanged loci are disjoint, the type-(II) incidence has no $\Fb_q$-points.

The type-(I) locus supports the order-four automorphisms whose square is inversion.  Its contribution after removal of the minimality defect is
\[
\alpha_5(q)q^{5M}-\beta_4(q)q^{4M}
=
\alpha_5(q)B^{5/12}-\beta_4(q)B^{1/3}.
\]
When $r$ is even, the four type-(II) fixed-pair loci support the order-three and order-six automorphisms.  Their combined contribution is
\[
4\epsilon_r\bigl(\alpha_6(q)q^{6M}-\beta_3(q)q^{3M}\bigr)
=
4\epsilon_r\alpha_6(q)B^{1/2}
-4\epsilon_r\beta_3(q)B^{1/4}.
\]

The additional nonsmooth equations in the quasi-elliptic-type stratum of Proposition~\ref{prop:nonsmooth} contribute the term $-2\beta_6(q)B^{1/2}$ in $2\Nc^w(K,B)$.  Thus the total coefficient of $B^{1/2}$ is $2\bigl(2\epsilon_r\alpha_6(q)-\beta_6(q)\bigr)$.  When $r$ is odd, this coefficient is the negative quantity $-2\beta_6(q)$.

As in characteristic three, the negative terms $-\beta_4(q)B^{1/3}$ and $-4\epsilon_r\beta_3(q)B^{1/4}$ are created by removing nonminimal models from positive ambient fixed-pair counts.

\subsection{Terms common to all characteristics}

For $M\geq1$ and $B=q^{12M}$, the calculations above, together with \cite[Theorem~9.7]{BPS}, show that in every characteristic the leading term
\[
2\left(\frac{q^9-1}{q^8-q^7}\right)B^{5/6}
\]
of $\Nc(K,B)$ is twice the corresponding leading term of $\Nc^w(K,B)$.  Likewise, the term $-2B^{1/6}=-2q^{2M}$ in $\Nc(K,B)$ is twice the corresponding weighted term and arises from subtracting equations whose generic fiber has a $K$-rational singular point.  In characteristics two and three, these common contributions, together with the characteristic-specific analyses above, account for every term in Theorem~\ref{thm:main}.

\medskip

    \section*{Acknowledgements}

    Warm thanks to Dori Bejleri, Andrea Di Lorenzo, Scott Mullane, Tristan Phillips and Matthew Satriano for numerous helpful discussions. The author was partially supported by the ARC grant DP210103397 and the Sydney Mathematical Research Institute.


    \bibliographystyle{alpha}
    \bibliography{main.bib}

    \vspace{+17pt}
    
    \noindent Jun--Yong Park \enspace -- \enspace \texttt{june.park@sydney.edu.au} \\
    \textsc{School of Mathematics and Statistics, University of Sydney, Australia} \\

\end{document}